\magnification1120\overfullrule=0pt

\input epsf.sty

\def\cA{{\cal A}}

\def\cE{{\cal E}}
\def\cF{{\cal F}}
\def\cH{{\cal H}}
\def\cI{{\cal I}}
\def\cL{{\cal L}}

\def\cO{{\cal O}}

\def\C{{\bf C}}

\def\P{{\bf P}}

\centerline{\bf Periodicities in Linear Fractional Recurrences:}

\centerline{\bf Degree growth of birational surface maps}
\bigskip
\centerline{Eric Bedford and Kyounghee Kim}
\bigskip
\bigskip\centerline {\bf \S0.  Introduction}  
\medskip\noindent Given complex numbers $\alpha_{0},\dots,\alpha_{p}$ and $\beta_{0},\dots,\beta_{p}$, we consider the recurrence relation

$$x_{n+p+1}= {\alpha_{0}+\alpha_{1}x_{n+1}+\cdots+\alpha_{p}x_{n+p} \over \beta_{0}+\beta_{1}x_{n+1}+\cdots+\beta_{p}x_{n+p}}.\eqno(0.1)$$
Thus a $p$-tuple $(x_{1},\dots,x_{p})$ generates an infinite sequence $(x_{n})$.  We consider two equivalent reformulations in terms of rational mappings:  we may consider the mapping $f:{\bf C}^{p}\to{\bf C}^{p}$ given by
$$f(x_{1},\dots,x_{p}) =\left(x_{2},\dots,x_{p},{\alpha_{0}+\alpha_{1}x_{1}+\cdots+\alpha_{p}x_{p} \over \beta_{0}+\beta_{1}x_{1}+\cdots+\beta_{p}x_{p}}\right).\eqno(0.2)$$
Or we may use the imbedding $(x_{1},\dots,x_{p})\mapsto[1:x_{1}:\dots:x_{p}]\in{\bf P}^{p}$ into projective space and consider the induced map $f:{\bf P}^{p}\to{\bf P}^{p}$ given by 
$$f_{\alpha,\beta}[x_{0}:x_{1}:\dots:x_{p}] = [x_{0}\beta\cdot x:x_{2}\beta\cdot x:\dots:x_{p}\beta\cdot x:x_{0}\alpha\cdot x],\eqno(0.3)$$
where we write $\alpha\cdot x=\alpha_{0}x_{0}+\cdots+\alpha_{p}x_{p}$.  

Here we will study the degree growth of the iterates $f^{k}=f\circ\cdots\circ f$ of $f$.  In particular, we are interested in the quantity 
$$\delta(\alpha,\beta):=\lim_{k\to\infty}\left({\rm degree}(f^{k}_{\alpha,\beta})\right)^{1/k}.$$
A natural question is: for what values of $\alpha$ and $\beta$ can (0.1) generate a periodic recurrence?  In other words, when does (0.1) generate a periodic sequence $(x_{n})$ for all choices of $x_{1},\dots,x_{p}$?  This is equivalent to asking when there is an $N$ such that $f^{N}_{\alpha,\beta}$ is the identity map.  Periodicities in recurrences of the form (0.1) have been studied in [L, KG, KoL, GL, CL].  The question of determining the parameter values $\alpha$ and $\beta$ for which $f_{\alpha,\beta}$ is periodic has been known for some time and is posed explicitly in [GKP] and [GL, p.\ 161].  Recent progress in this direction was obtained in [CL].   The connection with our work here is that if $\delta(\alpha,\beta)>1$, then the degrees of the iterates of $f_{\alpha,\beta}$ grow exponentially, and $f_{\alpha,\beta}$ is far from periodic.

In the case $p=1$, $f$ is a linear (fractional) map of ${\bf P}^{1}$.  The question of periodicity for $f$ is equivalent to determining when a $2\times 2$ matrix is a root of the identity.  In this paper we address these questions in the case $p=2$.  In fact, our principal efforts will be devoted to determining $\delta(\alpha,\beta)$ for all of the mappings in the family above.  In order to remove trivial cases, we will assume throughout this paper that
$$\eqalign{(\alpha_{0},\alpha_{1},\alpha_{2})&{\rm\ is\ not\ a\ multiple\ of\ }(\beta_{0},\beta_{1},\beta_{2}),\cr
  (\alpha_{1},\beta_{1})&\ne(0,0), \ \ (\alpha_{2},\beta_{2})\ne (0,0), {\rm\ and\ \ }\cr
&(\beta_{1},\beta_{2})\ne(0,0).\cr}\eqno(0.4)$$
Note that if the first condition in (0.4) is not satisfied, then the right hand side of (0.1) is constant.
If the left hand part of the second condition (0.4) is not satisfied, then $f$ does not depend on $x_{1}$ thus has rank 1, which cannot be periodic.  If the right hand part of the second condition (0.4) is not satisfied, then $f^{2}$ is essentially the 1-dimensional mapping $\zeta\mapsto{\alpha_{0}+\alpha_{1}\zeta\over \beta_{0}+\beta_{1}\zeta}$.  If the third condition in (0.4) is not satisfied, then $f$ is linear.  In this case, the periodicity of $f$ is a question of linear algebra.

Since we consider all parameters satisfying (0.4), we must treat a number of separate cases.  By $V_{n}$ we denote the variety of parameters $(\alpha,\beta)$ such that
$$\eqalign{&\ \ \beta_{2}=0,\ {\rm and\ \ }f^{n}_{\alpha,\beta}(q) = p,\cr
 {\rm\ where\ \ }&p=[\beta_{1}\alpha_{2}-\beta_{2}\alpha_{1}:-\beta_{0}\alpha_{2}+\alpha_{0}\beta_{2}:\alpha_{1}\beta_{0}-\alpha_{0}\beta_{1}],\cr
 {\rm\ and \ \ } &q=[\beta_{1}(\beta_{1}\alpha_{2}):\beta_{1}(\alpha_{1}\beta_{0}-\alpha_{0}\beta_{1}):\alpha_{1}(\beta_{1}\alpha_{2}-\alpha_{1}\beta_{2})].\cr}\eqno(0.5)$$
The following two numbers are of special importance here:
$$\phi\ (\sim 1.61803{\rm\ golden\ mean)\ is\ the\ largest\ root\ of\ }x^{2}-x-1\eqno(0.6)$$
$$\delta_{\star}\ (\sim 1.32472){\rm\ is\ the\ largest\ root\ of\ }x^{3}-x-1\eqno(0.7)$$
\proclaim Theorem  1.  If $(\alpha,\beta)\notin \bigcup_{n\ge0}V_{n}$, then $\phi\ge\delta(\alpha,\beta)\ge \delta_{\star}>1$.  For generic $(\alpha,\beta)$, the dynamic degree is $\delta(\alpha,\beta)=\phi$.

In particular, we see that $f_{\alpha,\beta}$ has exponential degree growth in all of these cases.  The remaining possibilities are:
\proclaim Theorem  2.  If $(\alpha,\beta)\in V_{n}$ for some $n\ge0$, then there is a complex manifold $X=X_{\alpha,\beta}$ obtained by blowing up ${\bf P}^{2}$ at finitely many points, and $f_{\alpha,\beta}$ induces a biholomorphic map $f_{\alpha,\beta}:X\to X$. Further:
\vskip0pt If {$n=0$},    $f_{\alpha,\beta}$ is periodic of period 6.
\vskip0pt If {$n=1$},   $f_{\alpha,\beta}$ is periodic of period 5.
\vskip0pt If {$n=2$},   $f_{\alpha,\beta}$ is periodic of period 8.
\vskip0pt If {$n=3$},    $f_{\alpha,\beta}$ is periodic of period 12.
\vskip0pt If {$n=4$},    $f_{\alpha,\beta}$ is periodic of period 18.
\vskip0pt If {$n=5$},    $f_{\alpha,\beta}$ is periodic of period 30.
\vskip0pt If {$n=6$},     the degree of $f^{n}_{\alpha,\beta}$ is asymptotically quadratic in $n$.
\vskip0pt If {$n\ge7$},  $f_{\alpha,\beta}$ has exponential degree growth rate  $\delta(\alpha,\beta)=\delta_{n}>1$, which is given by the largest root of the polynomial $x^{n+1}(x^{3}-x-1)+x^{3}+x^{2}-1$.  Further, $\delta_{n}$ increases to $\delta_{\star}$ as $n\to\infty$.

The family of maps
$$(x,y)\mapsto(y, {a+y\over x})$$
 has been studied by several authors (cf.\  [L, KoL, KLR, GBM, CL]).   Within this family, the case  $a=0$ corresponds to $V_0$, $a=1$ corresponds to $V_1$, and all the rest belong to the case $V_6$ (see \S6). 

In the cases $n\ge7$, the entropy of  $f_{\alpha,\beta}$ is equal to $\log\delta_{n}$ by Cantat [C].  The number $\delta_{\star}$ is known (see [BDGPS, Chap.\ 7]) to be the infimum of all Pisot numbers.   Diller and Favre [DF] showed that if $g$ is a birational surface map which is not birationally conjugate to a holomorphic automorphism, then $\delta(g)$ is a Pisot number.  So the maps $f$ in the cases $n\ge7$ have smaller degree growth than any such $g$.  Note that projective surfaces which have automorphisms of positive entropy are relatively rare:   Cantat [C] shows that, except for nonminimal rational surfaces (like $X$ in Theorem  2), the only possibilities are  complex tori, K3 surfaces, or Enriques surfaces.

\centerline{\epsfxsize=2.2in\epsfbox{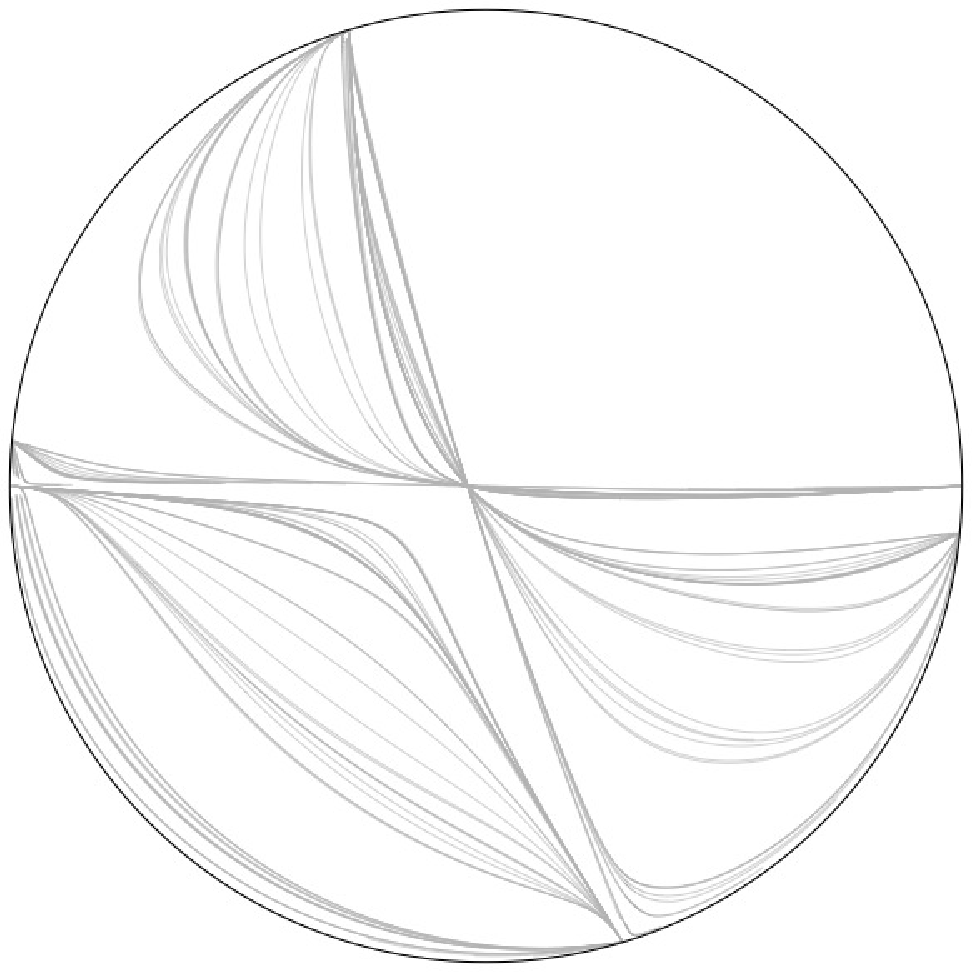}\hfil \epsfxsize=2.2in\epsfbox{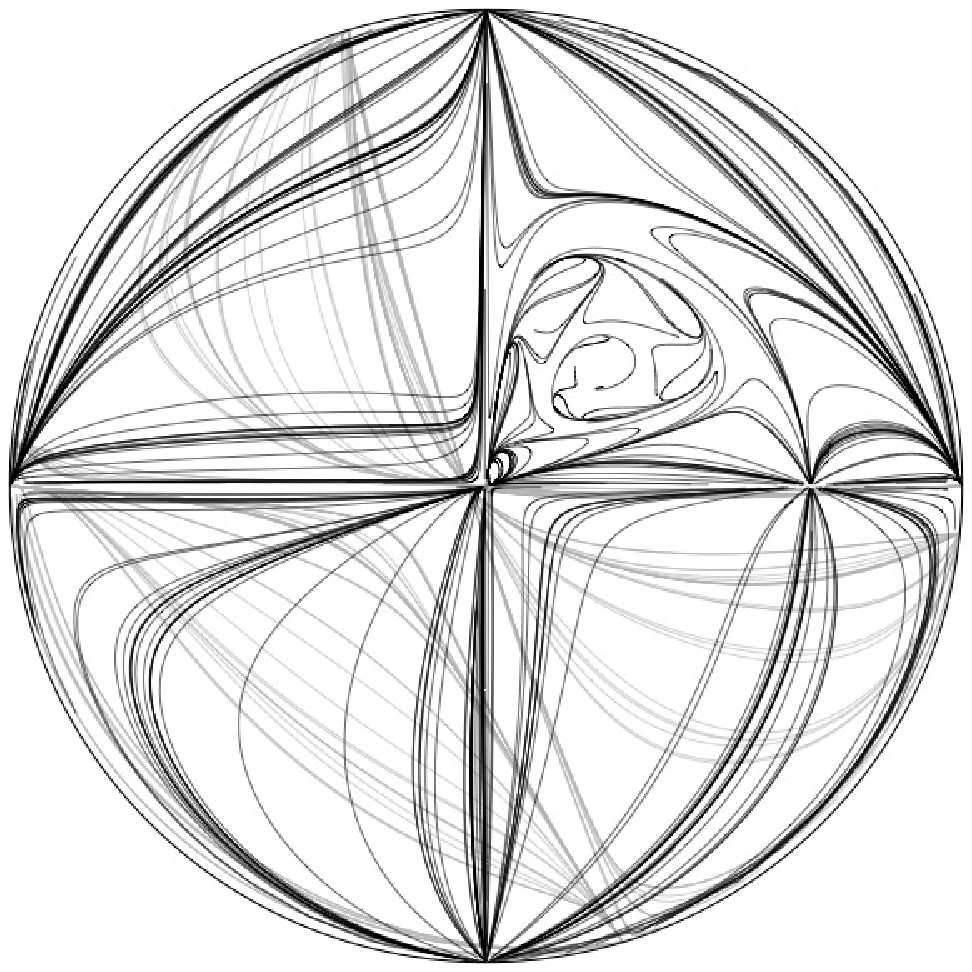}}
\centerline{Figure 0.1.  A map with (maximal) degree growth $\phi$.}
\medskip
Determining the dynamical degree for this family of mappings may be seen as a first step towards the dynamical study of these maps.  Figure 0.1 portrays the stable/unstable laminations of a mapping of maximal degree growth within the family $f_{\alpha,\beta}$. 
This paper is organized as follows.  In \S1 we give the general properties of the family $f_{\alpha,\beta}$.  In \S2 we show that $\delta(f_{\alpha,\beta})=\phi$ if $f_{\alpha,\beta}$ has only two exceptional curves.  Next we determine $\delta(f_{\alpha,\beta})$ in the (generic) case where it has three exceptional curves.  This determination, however, threatens to involve a large case-by-case analysis.  We avoid this by adopting a more general approach.  In \S3 we show how $\delta(f_{\alpha,\beta})$ may be derived from the set of numbers in open and closed orbit lists.  We do this by showing that results of [BK] may be extended from the ``elementary'' case to the general case.  We use this in \S4 to determine $\delta(\alpha,\beta)$ when the critical triangle is nondegenerate. In \S5 we handle the periodic cases in Theorem  2.  In \S6, we discuss parameter space and the varieties $V_{n}$ for $0\le n\le 6$.  We explain the computer pictures in the Appendix.
\medskip
We wish to thank Curt McMullen for helpful comments on this paper.

\bigskip
\centerline{\bf \S1.  Setting and Basic Properties}
\medskip
\noindent  In this section we review some basic properties of the map
$$f(x)=[x_{0}\beta\cdot x: x_{2}\beta\cdot x: x_{0}\alpha\cdot x],$$
which is the map (0.3) in the case $p=2$.  (We refer to [GBM] for a description of $f$ as a real map.)  The indeterminacy locus is
$$\eqalign{{\cal I} =& \{ x \in \P^2 : \, x_0 (\beta \cdot x) = x_2 (\beta \cdot x) = x_0 (\alpha \cdot x) = 0 \}\cr
=&\{e_{1},p_{0},p_{\gamma}\},} $$
where we set $p_{0}= [ 0:-\beta_2:\beta_1]$ and $p_{\gamma}=[\beta_{1}\alpha_{2}-\beta_{2}\alpha_{1}:-\beta_{0}\alpha_{2}+\alpha_{0}\beta_{2}:\alpha_{1}\beta_{0}-\alpha_{0}\beta_{1}]$.
Thus $f$ is holomophic on ${\bf P}^{2}-{\cal I}$, and its Jacobian  is  $2x_0(\beta \cdot x)[\beta_1 (\alpha \cdot x) - \alpha_1 ( \beta \cdot x)]$.  Let us set 
$$\gamma=(\beta_1 \alpha_0-\alpha_1\beta_0,0,\beta_1\alpha_2-\alpha_1\beta_2) \in \C^3$$
and note that the Jacobian vanishes on the curves
$$ \Sigma_0 = \{ x_0 = 0\}, \quad  \Sigma_\beta = \{ \beta \cdot x = 0 \},\quad {\rm and} \quad \Sigma_\gamma =\{\gamma \cdot x = 0\}.$$
These curves are exceptional in the sense that they are mapped to points:
$$f(\Sigma_{0}-{\cal I})=e_{1}:=[0:1:0],\ \ \ f(\Sigma_{\beta}-{\cal I})=e_{2}:=[0:0:1],\ \ f(\Sigma_{\gamma}-{\cal I})=q,\eqno(1.1)$$
where $q$ is defined in (0.5).  We write the set of exceptional curves as ${\cal E}(f)=\{\Sigma_{0},\Sigma_{\beta},\Sigma_{\gamma}\}$.

\proclaim Lemma 1.1.  $$ f(\P^2 - \Sigma_0 \cup\Sigma_\beta ) \cap \Sigma_0 = \emptyset.$$
Further, if $\beta_2\ne0$,
$$ f(\P^2 - \Sigma_0 \cup \{p_\gamma\}) \cap \{p_0\} = \emptyset.$$

\noindent{\it Proof.}  
In $\P^2 - \cE(f) \cup \cI(f)$, $f$ is holomorphic. It follows that for $[x_0:x_1:x_2] \in \P^2 - \cE(f) \cup \cI(f)$, $f([x_0:x_1:x_2]) \notin \Sigma_0$ since $x_0(\beta \cdot x) \ne 0$. If $\beta_1 =0$ or $\beta_1\alpha_2 - \alpha_1\beta_2 = 0$ then either $\Sigma_\gamma = \Sigma_\beta$ or $\Sigma_\gamma = \Sigma_0$. If both $\beta_1$ and $\beta_1\alpha_2 - \alpha_1\beta_2 $ are non-zero,  we have $f(\Sigma_\gamma)=q \notin \Sigma_0$.  In case $\beta_2 \ne 0$, for $[x_0:x_1:x_2] \in \Sigma_\beta$, we have seen that $f([x_0:x_1:x_2])= e_2 \ne p_0$, which completes the proof.
\medskip
The inverse of $f$ is given by the map
$$f^{-1}(x)=[x_{0}B\cdot x:x_{0}A\cdot x-\beta_{2}x_{1}x_{2}:x_{1} B\cdot x],$$
where we set $A=(\alpha_{0},\alpha_{2},-\beta_{0})$ and $B=(-\alpha_{1},0,\beta_{1})$.  In the special case $\beta_{2}=0$, the form of $f^{-1}$ is similar to that of $f$.  The indeterminacy locus ${\cal I}(f^{-1})=\{e_{1},e_{2},q\}$ consists of the three points which are the  $f$-images of the exceptional lines for $f$.  The Jacobian of $f^{-1}$ is 
$$-2x_{0} B\cdot x(\alpha_{1}\beta_{0}x_{0}-\alpha_{0}\beta_{1}x_{0}-\alpha_{2}\beta_{1}x_{1}+\alpha_{1}\beta_{2}x_{1}).$$
 Let us set  $C=(\alpha_{1}\beta_{0}-\alpha_{0}\beta_{1},\alpha_{1}\beta_{2}-\alpha_{2}\beta_{1},0)$, and $\Sigma_{B}=\{x\cdot B=0\}$,  $\Sigma_{C}=\{x\cdot C=0\}$. 
In fact, ${\cal E}(f^{-1})=\{\Sigma_{0},\Sigma_{B},\Sigma_{C}\}$, and  $f^{-1}$ acts as:  $\Sigma_{0}\mapsto p_{0}$, $\Sigma_{B}\mapsto e_{1}$, and $\Sigma_{C}\mapsto p_{\gamma}$.

To understand the behavior of $f$ at ${\cal I}$, we define the cluster set $Cl_{f}(a)$ of a point  $a \in \P^2$ by $$Cl_f(a) = \{ x \in \P^2 : \, x = \lim_{a' \to a} f(a'), \,\, a' \in \P^2 -\cI(f) \}.$$
In general, a cluster set is connected and compact.  In our case, we see that the cluster set is a single point when $a\notin{\cal I}$, i.e., when $f$ is holomorpic.  And the cluster sets of the points of indeterminacy are found by applying $f^{-1}$: i.e.,  $e_{1}\mapsto Cl_{f}(e_{1})=\Sigma_{B}$, $p_{0}\mapsto Cl_{f}(p_{0})=\Sigma_{0}$, and $p_{\gamma}\mapsto Cl_{f}(p_{\gamma})=\Sigma_{C}$.  Thus $f$ acts as in Figure 1.1: the lines on the left hand triangle are exceptional and are mapped to the vertices of the right hand triangle, and the vertices of the left hand triangle are blown up to the sides of the right hand triangle.

Let 
$$\pi:Y\to{\bf P}^{2}\eqno(1.2)$$ 
be the complex manifold obtained by blowing up ${\bf P}^{2}$ at $e_{1}$.  We will discuss the induced birational map $f_{Y}:Y\to Y$.  We let $E_{1}:=\pi^{-1}e_{1}$ denote the exceptional blow-up fiber.  The projection gives a biholomorphic map $\pi:Y-E_{1}\to {\bf P}^{2}-e_{1}$.  For a complex curve $\Gamma\subset {\bf P}^{2}$, we use the notation $\Gamma\subset Y$ to denote the {\it strict transform} of $\Gamma$ in $Y$.  Namely, $\Gamma$ denotes the closure of $\pi^{-1}(\Gamma-e_{1})$ inside $Y$.  Thus $\Gamma$ is a proper subset of $\pi^{-1}\Gamma=\Gamma\cup E_{1}$.

We identify $E_{1}$ with ${\bf P}^{1}$ in the following way.  For $[\xi_{0}:\xi_{2}]\in{\bf P}^{1}$, we associate the point
$$[\xi_{0}:\xi_{2}]_{E_{1}}:=\lim_{t\to0}\pi^{-1}[t\xi_{0}:1:t\xi_{2}]\in E_{1}.$$
We may now determine the map $f_{Y}$ on $\Sigma_{0}$.  For $x=[0:x_{1}:x_{2}]=\lim_{t\to0}[t:x_{1}:x_{2}]\in\Sigma_{0}$, we assign $f_{Y}x:=\lim_{t\to0}f[t:x_{1}:x_{2}]\in Y$.  That is, 
$f[t:x_{1}:x_{2}]=[t\beta\cdot x:x_{2}\beta\cdot x:t\alpha\cdot x]$, and so taking the limit as $t\to0$, we obtain
$$f_{Y}[0:x_{1}:x_{2}]=[\beta\cdot x:\alpha\cdot x]_{E_{1}}.\eqno(1.3)$$

\epsfxsize=4in
\centerline{\epsfbox{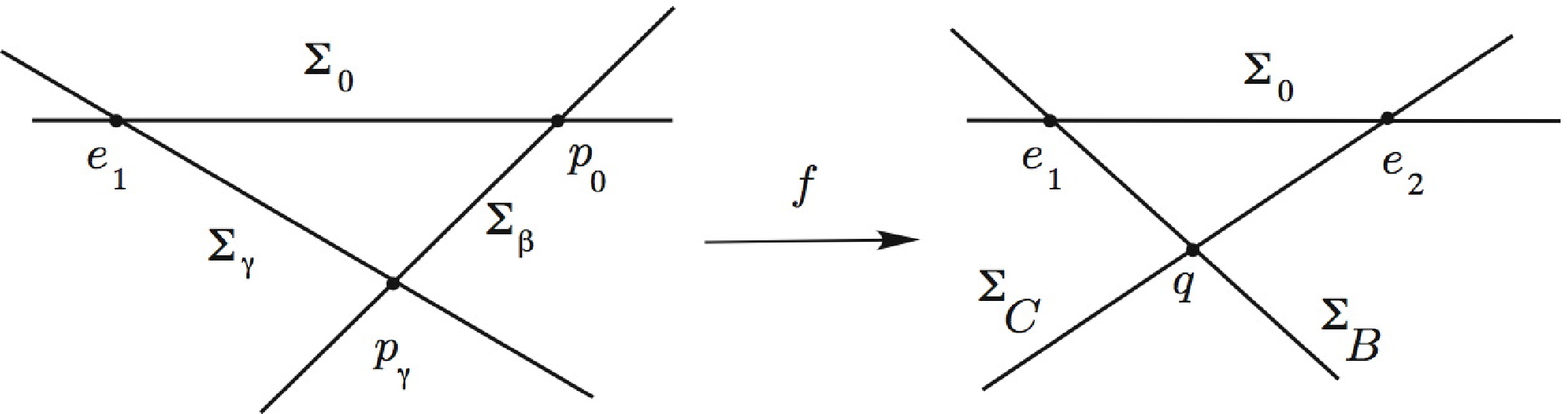}}
$$ p_{0}\mapsto\Sigma_{0}\mapsto e_{1}\mapsto\Sigma_{B}\qquad
\Sigma_{\beta}\mapsto e_{2}\qquad
\Sigma_{\gamma}\mapsto q \qquad p_{\gamma}\mapsto\Sigma_{C}$$
\centerline{Figure 1.1.  Blowing-up/blowing-down behavior of $f$. }
\medskip

Now we make a similar computation for a point $[\xi_{0}:\xi_{2}]_{E_{1}}$ in the fiber $E_{1}$ over the point of indeterminacy $e_{1}$.  We set $x=[t\xi_{0}:1:t\xi_{2}]$ so that
$$fx = [t\xi_{0}\beta\cdot x:t\xi_{2}\beta\cdot x:t\xi_{0}\alpha\cdot x].$$
Taking the limit as $t\to0$, we find 
$$f_{Y}([\xi_{0}:\xi_{2}]_{E_{1}})= [\xi_{0}\beta_{1}:\xi_{2}\beta_{1}:\xi_{0}\alpha_{1}]\in\Sigma_{A}.$$
Thus we have:

\proclaim Lemma 1.2.  The map $f_{Y}$ has the properties:
 \item{(i)}$f_{Y}$ is a local diffeomorphism at points of $\Sigma_{0}$ if and only if $\beta_{1}\alpha_{2}-\alpha_{1}\beta_{2}\ne0$.
\item{(ii)}  $f_{Y}$ is a local diffeomorphism at points of $E_{1}$ if and only if $\beta_{1}\ne0$.

\medskip
\centerline{\bf \S2.  Degenerate Critical Triangle}
\medskip
\noindent We will refer to the set $\{\Sigma_{0},\Sigma_{\beta},\Sigma_{\gamma}\}$ of  exceptional curves as the {\it critical triangle}; we say that the critical triangle is {\it nondegenerate} if these three curves are distinct.  Since $(\beta_{1},\beta_{2})\ne(0,0)$, we have $\Sigma_{0}\ne\Sigma_{\beta}$.  Thus there are only two possibilities for a degenerate triangle.  The first of these is the case $\Sigma_{\gamma}=\Sigma_{\beta}$, which occurs when $\beta_{1}=0$.  The second is $\Sigma_{\gamma}=\Sigma_{0}$, which occurs when $\beta_{1}\alpha_{2}-\alpha_{1}\beta_{2}=0$.  (And since $\Sigma_{0}\ne\Sigma_{\beta}$ we have $\beta_{1}\ne0$ in this case.)  We will show that $\delta(\alpha,\beta)=\phi$ when the critical triangle is degenerate.  This is different from the general case (and easier), and we treat it in this section.

In order to determine the degree growth rate of $f$, we will consider the induced pullback $f^{*}$ on $H^{1,1}$.  We will be working on compact, complex surfaces $X$ for which $H^{1,1}(X)$ is generated by the classes of divisors.  If $[D]$ is the divisor of a curve $D\subset X$, then we define $f^{*}[D]$ to be the class of the divisor $f^{-1}D$.  We say that $f$ is {\it 1-regular} if $(f^{n})^{*}=(f^{*})^{n}$ for all $n\ge0$.  Fornaess and Sibony showed in [FS] that if
$${\rm\ for\ every\ exceptional\ curve\ }C{\rm\ and\ all\ }n\ge0, f^{n}C\notin{\cal I}\eqno(2.1)$$
then $f$ is 1-regular.  We will use this criterion in the following:  
\proclaim Proposition 2.1.  If the critical triangle is degenerate, then the map $f_{Y}:Y\to Y$ is 1-regular.

\noindent{\it Proof.}  We treat the two possibilities separately.  The first case is $\Sigma_{\gamma}=\Sigma_{\beta}$; see Figure 2.1.  In this case $f$ has two exceptional lines $\Sigma_{0}$ and $\Sigma_{\beta}$ and two points of indeterminacy ${\cal I}=\{e_{1},p_{\lambda}\}$.  After we blow up $e_{1}$ to obtain $Y$, the line $\Sigma_{0}$ is no longer exceptional.  (Our drawing convention in this and subsequent Figures is that exceptional curves are thick, and points of indeterminacy are circled.)  By (1.3), we see that $f_{Y}$ maps $E_{1}$ to $e_{2}=q$, and thus the exceptional set becomes ${\cal E}(f_{Y})=\{E_{1},\Sigma_{\beta}=\Sigma_{\gamma}\}$.  Now in order to check condition (2.1), we need to follow the orbit of $e_{2}$.  By (1.3) we see that $e_{2}$ is part of a 2-cycle $\{e_{2},[\beta_{2}:\alpha_{2}]_{E_{1}}\}$.  On the other hand, the points of indeterminacy for $f_{Y}$ are $p_{\gamma}$ and $[0:1]_{E_{1}}=E_{1}\cap\Sigma_{0}$.  Since $\beta_{1}=0$ in this case, we have $\beta_{2}\ne0$, so (2.1) holds.

\epsfxsize=4in
\centerline{\epsfbox{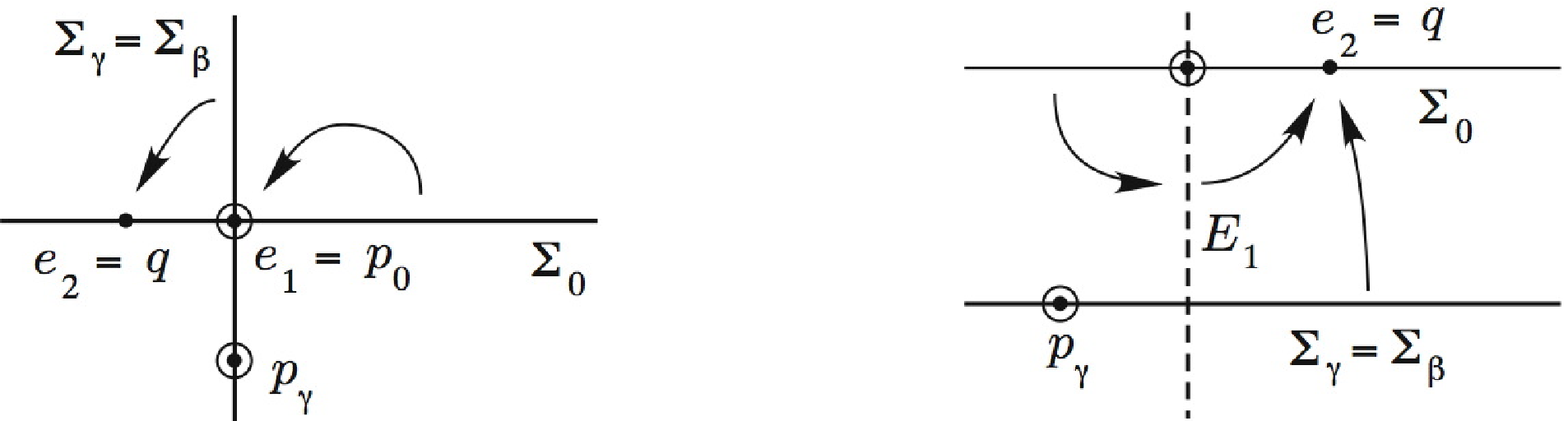}}
\centerline{Figure 2.1.  The case $\Sigma_{\beta}=\Sigma_{\gamma}$.}
\medskip

The second case is $\Sigma_{\gamma}=\Sigma_{0}$.  Again, ${\cal I}
=\{e_{1},p_{\gamma}\}$, but ${\cal E}(f)=\{\Sigma_{0},\Sigma_{\beta}\}$, and the arrangement of exceptional curves and points of indeterminacy are as in Figure 2.2.  In this case, we have $\beta_{1}\ne0$, so by Lemma 1.2, we hve ${\cal I}(f_{y})=\{p_{0}=p_{\gamma}\}$ and ${\cal E}(f_{Y})=\{\Sigma_{\beta}\}$.  As before, we need to track the orbit of $e_{2}$.  But by Lemma 1.1, we see that we can never have $f^{j}e_{2}=p_{0}$ for $j\ge1$.  Thus (2.1) holds in this case, too, and the proof is complete.

\epsfxsize=4in
\centerline{\epsfbox{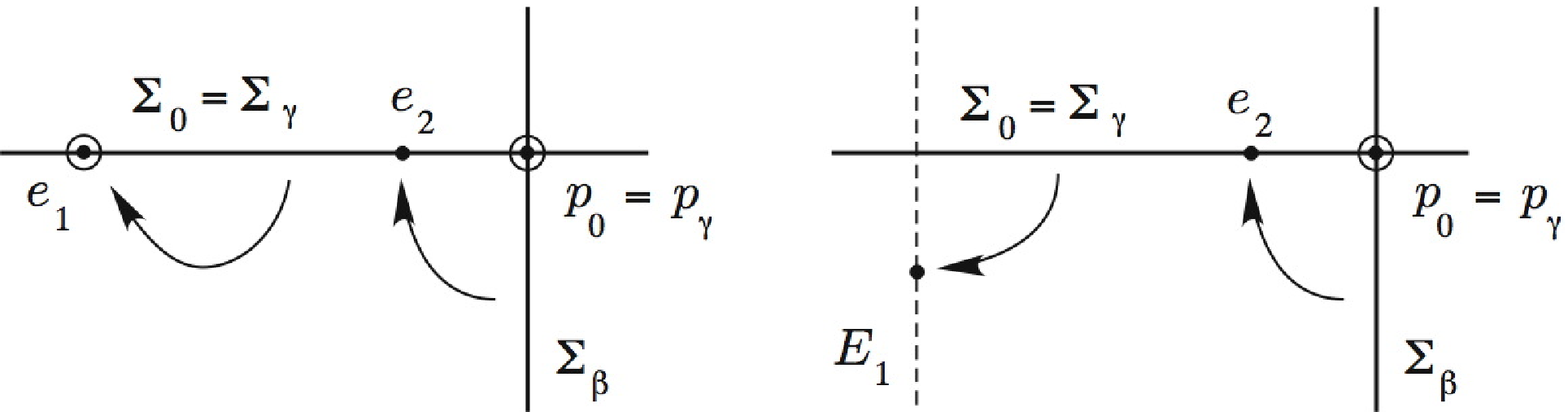}}
\centerline{Figure 2.2.  The case $\Sigma_{0}=\Sigma_{\gamma}$.}
\medskip

Now let us determine $f^*_Y$.  The cohomology group $H^{1,1}({\bf P}^{2};{\bf Z})$ is one-dimensional and is generated by the class of a complex line.  We denote this generator by $L$.  Let $L_{Y}:=\pi^{*}L\in H^{1,1}(Y;{\bf Z})$ be the class induced by the map (1.2).  It follows that $\{L_{Y},E_{1}\}$ is a basis for  $H^{1,1}(Y;{\bf Z})$.  Now $\Sigma_{0}=L\in H^{1,1}({\bf P}^{2};{\bf Z})$.  Pulling this back by $\pi$, we have
$$L_{Y}=\pi^{*}\Sigma_{0}=\Sigma_{0}+E_{1}.$$
Now $f_{Y}^{*}$ acts by taking pre-images:
$$f_{Y}^{*}E_{1}=[f^{-1}E_{1}] = \Sigma_{0}=L_{Y}-E_{1},$$
where the last equality follows from the equation above.

Now $e_{1}$ is indeterminate, and $fe_{1}=\Sigma_{A}$.  Since $\Sigma_{A}$ intersects any line $L$, it follows that $e_{1}\in f^{-1}L$.  Thus
$$\pi^{*}[f^{-1}L]=[f^{-1}L] + E_{1}\in H^{1,1}(Y;{\bf Z}).$$
On the other hand, $f^{-1}L=2L\in H^{1,1}({\bf P}^{2};{\bf Z})$.  Thus
$$\pi^{*}[f^{-1}L]=\pi^{*}2L=2L_{Y}.$$
Putting these last two equations together, we have $f_{Y}^{*}L_{Y}=2L_{Y}-E_{1}$.  Thus
$$f^{*}_{Y}=\pmatrix{2&1\cr -1&-1\cr},$$
which is a matrix with spectral radius equal to $\phi$.  This yields the following:
\proclaim Proposition 2.2.  If the critical triangle is degenerate, then $\delta(\alpha,\beta)=\phi$.

\centerline{\bf \S3.  Regularization and Degree Growth}
\medskip\noindent 
In this Section we discuss a different, but more general, family of maps.  By $J:{\bf P}^2\to{\bf P}^2$ we denote the involution
$$J[x_0:x_1:x_2]=[x_0^{-1}:x_1^{-1}:x_2^{-1}]=[x_1x_2:x_0x_2:x_0x_1].$$
For an invertible linear map $L$ of $\P^2$ we consider the map $f:=L\circ J.$ The exceptional curves are 
$\cE=\{{\bf \Sigma}_0,{\bf \Sigma}_1,{\bf \Sigma}_2\},$ where ${\bf \Sigma}_j : =\{x_j =0\}, j=0,1,2,$ and the points of indeterminacy are $\cI=\{\epsilon_0,\epsilon_1,\epsilon_2\},$ where $\epsilon_i={\bf \Sigma}_j \cap {\bf \Sigma}_k,$ with $\{i,j,k\}=\{0,1,2\}.$ We define ${\bf a}_j :=f({\bf\Sigma}_j-\cI) = L\epsilon_j$ for $j=0,1,2.$

For $p\in \P^2$ we define the {\it orbit} $\cO(p)$ as follows. If $p \in \cE \cup \cI$, then $\cO(p) = \{p\}.$ If there exists an $N\geq 1$ such that $f^jp \notin \cE\cup\cI$ for $0 \leq j \leq N-1$ and $f^Np \in \cE\cup\cI,$ then we set $\cO(p) = \{p,fp,\dots,f^Np\}.$ Otherwise we have $f^jp\notin\cE\cup\cI$ for all $j \geq0,$ and we set $\cO(p) = \{p, fp,f^2p,\dots \}.$ We say the orbit $\cO(p)$ is {\it singular} if it is finite; otherwise, it is non-singular. We say an orbit $\cO(p)$ is {\it elementary} if it is either non-singular, or if it ends at a point of indeterminacy.  In other words, a non-elementary orbit ends in a point of ${\cal E}-{\cal I}$. 

\proclaim Lemma 3.1.  If $f$ has at least one singular orbit, then it has a singular orbit that is elementary.

\noindent{\it Proof.}  Suppose for all $i \in S_0$, $\cO_i$ is non-singular. If follows that every orbit $\cO_j, j \notin S_0$ ends at a point in ${\bf \Sigma}_i, i \in S_0$. Since all $\cO_i, i \in S_0$ are non-singular, ${\bf\Sigma}_j, j \notin S_0$ cannot end at a point of indeterminacy. This means that $f$ is $1$-regular.
\bigskip
Henceforth, we will assume that $f$ has singular orbits.  Let us write $\cO_i=\cO({\bf a}_i)=\cO(f({\bf \Sigma}_i-\cI))$ for the orbit of an exceptional curve.  We set 
$$S=\{i\in\{0,1,2\}: \cO_i{\rm\ is\ singular}\},$$
and
$$S_0=\{i\in\{0,1,2\}:\cO_i{\rm\ is\ singular\ and\ elementary}\}.$$
 Let $\cO_{S_0}=\bigcup_{i\in S_0}\cO_i$.  We write $X_0={\bf P}^2$, and let $\pi:X_1\to X_0$ be the complex manifold obtained by blowing up the points of $\cO_{S_0}$.  We let $f_1:X_1\to X_1$ denote the induced birational mapping. By Lemma 1.2, we see that the curves ${\bf \Sigma}_i$, $i\in S_0$, are not exceptional for $f_1$, and the blowing up operation constructed no new points of indeterminacy for $f_1$.  Thus the exceptional curves for $f_1$ are ${\bf \Sigma}_i$ for $i\notin S_0$.   If $S_0$ is a proper subset of $S$, then for $i\in S-S_0$ we redefine $\cO_i$ to be the $f_1$-orbit of ${\bf a}_i$ inside $X_1$.   Let us define $S_1=\{i\in S-S_0: \cO_i{\rm\ is\ elementary}\}$.  We may apply Lemma 3.1 to conclude that if $S-S_0\ne\emptyset$, then $S_1\ne\emptyset$.  As before, we may define $\cO_{S_1}=\bigcup_{i\in S_1}\cO_i$, and we construct the complex manifold $\pi:X_2\to X_1$ by blowing up all the points of $\cO_{S_1}$.  Doing this, we reach the situation where every singular orbit $\cO_i$ has the property that it is elementary in some $X_j$, and thus it has the form $\cO_i=\{{\bf a}_i,\dots,\epsilon_{\tau(i)}\}$ for some $\tau(i)\in\{0,1,2\}$.  

Next we organize the singular orbits $\cO_i$  into lists, as follows (modulo permutation of the indices $\{0,1,2\}$).  If there is only one singular orbit, we have the list $\cL=\{\cO_i=\{{\bf a}_i, \dots ,\epsilon_{\tau(i)}\}\}.$ If $\tau(i)=i,$ we say that $\cL$ is a {\it closed list}; otherwise it is an {\it open list}. If there are two singular orbits, we can have two closed lists:
$$\cL_1=\{\cO_0=\{{\bf a}_0,\dots,\epsilon_0\}\}, \ \cL_2 = \{\cO_1= \{{\bf a}_1, \dots,\epsilon_1\}\}$$
or a closed list and an open list:
$$\cL_1=\{\cO_0=\{{\bf a}_0,\dots,\epsilon_0\}\}, \ \cL_2 = \{\cO_1= \{{\bf a}_1, \dots,\epsilon_2\}\}.$$
We cannot have two open lists since there are only $3$ orbits $\cO_i$. We can also have a single list :
$$\cL=\{\cO_0=\{{\bf a}_0, \dots,\epsilon_1\},\cO_1=\{{\bf a}_1,\dots,\epsilon_{\tau(1)}\}\},$$
which is a closed list if $\tau(1)=0$ and an open list otherwise. If there are three singular orbits, then the possibilities are 
$$\cL=\{\cO_0=\{{\bf a}_0,\dots,\epsilon_1\},\cO_1=\{{\bf a}_1,\dots,\epsilon_2\}, \cO_2=\{{\bf a}_2, \dots,\epsilon_0\}\},$$
$$\cL_1=\{\cO_0=\{{\bf a}_0,\dots,\epsilon_0\}\},\ \cL_2=\{\cO_1=\{{\bf a}_1,\dots,\epsilon_2\}, \cO_2=\{{\bf a}_2, \dots,\epsilon_1\}\},$$
or
$$\cL_1=\{\cO_0=\{{\bf a}_0,\dots,\epsilon_0\}\}, \cL_2=\{\cO_1=\{{\bf a}_1,\dots,\epsilon_1\}\}, \cL_3=\{\cO_2=\{{\bf a}_2,\dots,\epsilon_2\}\},$$
where all the lists are closed.

For an orbit $\cO_i$, we let $n_i= |\cO_i|$ denote its length, and for an orbit list $\cL= \{\cO_a, \dots, \cO_{a+\mu}\},$ we denote the set of orbit lengths by $|\cL|= \{n_a, \dots, n_{a+\mu}\}.$ We set $\#\cL^c  = \{ |\cL_j| : \cL_j {\rm \ is \ closed} \}$ and  $\#\cL^o  = \{ |\cL_j| : \cL_j {\rm \ is \ open} \}.$ The set $\#\cL^c$ and $\#\cL^o$ determine $\delta(f)$, as is shown in the following:

\proclaim Theorem 3.2. If $f=L\circ J$, then the dynamic degree $\delta(f)$ is the largest real zero of the polynomial
$$\chi(x) = (x-2) \prod_{\cL \in \cL^c\cup \cL^o} T_{\cL}(x)+ (x-1) \sum_{\cL \in \cL^c\cup \cL^o}  S_{\cL}(x) \prod_{\cL' \ne \cL}  T_{\cL'}(x). \eqno{(3.1)}$$
Here $\cL$ runs over all orbit lists. For each orbit list $\cL$, we let $N$ denote the sum of all the length of the orbits in $\cL$. If  $\cL$ is closed $T_{\cL}(x)= x^{N}-1$, and if $\cL$ is open $T_{\cL}(x)= x^{N}$. The polynomial $S_{\cL}$ is defined by
$$ \eqalign{ S_{\cL} (x) \ &= 1 \cr &= x^{n_1}+x^{n_2}+2 \cr & =  x^{n_1}+x^{n_2}+1\cr &=\sum^3_{i=1}[x^{N-n_i}+x^{n_i}]+3\cr &=\sum^3_{i=1}x^{N-n_i}+\sum_{i\ne 2}x^{n_i}+1\cr} \quad \eqalign{ &{\rm if\ } |\cL|=\{n_1\} \cr   & {\rm if \ } \cL {\rm \ is \ closed \ and \ }  |\cL|=\{n_1, n_2\} \phantom{x^N}\cr   & {\rm if \ } \cL {\rm \ is \ open\ and \ }  |\cL|=\{n_1, n_2\} \phantom{x^N} \cr   & {\rm if \ } \cL {\rm \ is \ closed \ and \ }  |\cL|=\{n_1, n_2,n_3\}  \phantom{\sum^3_{i=1}} \cr  & {\rm if \ } \cL {\rm \ is \ open \ and \ }  |\cL|=\{n_1, n_2,n_3\}. \phantom{\sum^3_{i=1}}\cr }$$

The rest of this section is devoted to the proof of Theorem $3.2.$ We start by considering the case where $f$ is elementary. In this case we have $S=S_0$.  We set $X:=X_1$.  It follows from (2.1) that $f_X:X\to X$ is 1-regular, and thus $\delta(f)$ is spectral radius of $f_X^*$.  The computation given in the Appendix of [BK] then shows that $(3.1)$ is the characteristic polynomial of $f_X^*$.

For $p\in \cO_S-\cI$ we let $\cF_p=\pi^{-1}p$ denote the exceptional fiber over $p$.  If $\epsilon_i\in\cO_S\cap\cI$, we let $E_i$ denote the exceptional fiber over $\epsilon_i$.   We will feel free to identify curves with the classes they generate in $H^{1,1}(X)$.  Let $H \in H^{1,1}(\P^2)$ denote the class of a line, and let $H_X = \pi^*H$ denote the induced class in $H^{1,1}(X)$. For $i \in S$, we have
$$ {\bf \Sigma}_i \to {\bf a}_i \to \cdots \to f^{n_i-1}{\bf a}_i = f^{n_i}({\bf \Sigma}_i - \cI) = \epsilon_{\tau(i)}$$
for some $\tau(i) \in \{0,1,2\}$. At each points $f^j{\bf a}_i, \ 0 \leq j \leq n_i-1$, $f$ is locally biholomorphic, so $f_X$ induces a biholomorphic map 
$$f_X : \cF_{f^j {\bf a}_i} \to \cF_{f^{j+1} {\bf a}_i} \quad \ 0 \leq j \leq n_i-2, {\rm \  and}$$
$$f_X : \cF_{f^{n_i-1} {\bf a}_i} \to E_{\tau(i)}.$$
It follows that
$$\eqalign{ & f_X^* \cF_{f^{j+1} {\bf a}_i} = \cF_{f^j {\bf a}_i} \quad {\rm for \ }\  0 \leq j \leq n_i-2,\ i \in S \cr & f_X^* E_{\tau(i)} = \cF_{f^{n_i-1} {\bf a}_i}\cr} \eqno{(3.2)}$$
and 
$$f_X^*\cF_{{\bf a}_i} = \{{\bf \Sigma}_i\} \quad {\rm for \ } i \in S \eqno{(3.3)}$$
where $\{{\bf \Sigma}_i\}$ is the induced class by ${\bf \Sigma}_i$ in $H^{1,1}(X)$. Let $\Omega = \cI \cap \{ \epsilon_{\tau(i)}=f^{n_i-1}{\bf a}_i, \ i \in S\}$, the set of blow-up centers which belongs to $\cI$. Let us denote $\cA$ the set of indices $i$ such that $\cO_i$ is singular orbit and is the first orbit in an open orbit list. For each $i$, ${\bf \Sigma}_i$ contains blow-up centers in the set $\Omega- \{\epsilon_i\}$. Notice that if $i \in \cA$, $\epsilon_i \notin \Omega$, otherwise $\epsilon_i \in \Omega$. Using the identigy $\pi^*\{{\bf \Sigma}_i\} = \{ \pi^{-1} {\bf \Sigma}_i\}$, we have 
$$\eqalign{ & \{{\bf \Sigma}_i\} = H_X - E_\Omega + E_i \quad i \notin \cA\cr & \{{\bf \Sigma}_i\} = H_X - E_\Omega  \quad i \in \cA\cr}\eqno{(3.4)}$$ where $E_\Omega := \sum_{\epsilon_t \in \Omega} E_t $. A generic hyperplane $\cH$ in $\P^2$ does not contain any blow-up centers and may be considered to be subset of $X$. Let us restrict the map to $X-\cI$. A generic hyperplane $\cH$ does intersect with any line in $\P^2$. It follows that $\epsilon_i \in f^{-1}_X \cH$, $ i \in \Omega$ and we have
$$2H_X = \pi^*(f^*H)= \pi^*\{f^{-1}\cH\} = f^*_X+ E_\Omega.$$
Therefore under $f^*_X$, we have
$$ f_X^*H_X= 2H_X-E_\Omega. \eqno{(3.5)}$$

Now let us suppose that $f$ is not elementary.  Let $S=S_0\cup S_1\cup S_2$ and the manifolds $\pi_{i+1}:X_{i+1}\to X_i$  be as above.  Let us set $X:=X_3$.  By (2.1) again, the induced map $f_X:X\to X$ is 1-regular.  If $p\in X_i$ is a center of blow-up, we let $\cF_p$ denote the exceptional fiber inside $X_{i+1}$, and we use the same notation for the divisor in $X$ given by the strict transform of $\cF_p$; in particular, $\cF_p$ is irreducible.  Thus for $i \in S_k$, $f_X$ induces a biholomorphic map 
$$f_X : \cF_{f_k^j {\bf a}_i} \to \cF_{f_k^{j+1} {\bf a}_i} \quad \ 0 \leq j \leq n_i-1.$$
It follows that 
$$\eqalign{ & f_X^* \cF_{f_k^{j+1} {\bf a}_i} = \cF_{f_k^j {\bf a}_i} \quad {\rm for \ }\  0 \leq j \leq n_i-1,\ i \in S_k \cr & f_X^* \cF_{{\bf a}_i} = \{{\bf\Sigma}_i\}  {\rm \ for\ }i \in S} \eqno{(3.6)}$$
where $\{{\bf \Sigma}_i\}$ is the induced class by ${\bf \Sigma}_i$ in $H^{1,1}(X)$. Let $\Omega_i:=\{p\in \cO_S: \pi(p)= \epsilon_i\}$ the set of blow-up centers whose image of $\pi$ is $\epsilon_i\in\cI$ and let $\Omega:=\bigcup_{i \in S}\Omega_i.$ For each $i\in S$, we denote $\Xi_i:=\{p\in \cO_S: \pi(p)\in {\bf\Sigma}_i - \cI\}$ the set of blow-up centers which belongs to exceptional line ${\bf\Sigma}_i - \cI$ and we let $\Xi = \bigcup_{i \in S}\Xi_i.$ We also use the notation $\cA$ for the set of indices $i$ such that $\cO_i$ is singular orbit and is the first orbit in an open orbit list. For each $i\in S$, ${\bf \Sigma}_i$ contains blow-up centers in the set $(\Omega-\Omega_i)\cup \Xi_i$. with $H_X = \pi^*H$ the induced class in $H^{1,1}(X)$, we have 
$$\eqalign{ & \{{\bf \Sigma}_i\} = H_X - E_\Omega + E_i-\cF_{\Xi_i} \quad i \notin \cA\cr & \{{\bf \Sigma}_i\} = H_X - E_\Omega -\cF_{\Xi_i}  \quad i \in \cA\cr}\eqno{(3.7)}$$
where $E_\Omega=\sum_{p\in \Omega} \cF_p$, $E_i=\sum_{p\in \Omega_i} \cF_p$, and $\cF_{\Xi_i}=\sum_{p\in \Xi_i} \cF_p.$ We also have
$$2H_X = \pi^*(f^*H)=  f^*_X+ E_\Omega.\eqno{(3.8)}$$

To finish the proof, let us suppose that $g$ is an elementary map, and $f$ is not elementary, but both have the same orbit list structure given by $\#\cL^c, \#\cL^o$.  We have shown that $g^*$ is represented by the transformation (3.2--5), and $f^*$ is represented by the transformation (3.6--8).  To finish the proof, we show that these two linear transformations have the same characteristic polynomials.   We illustrate this computation with an example which appears later in the paper. (The matrix computation for the other cases are similar.)  We consider the case where the list structures of $f$ and $g$ are both given by 
$$\#\cL^o= \emptyset, \quad \#\cL^c= \{ \{1,6\}\}.$$
We may also assume that $1 \in S_0$, $2\in S_1$ and $$\cO_1= \{ {\bf a}_1=\epsilon_2\}, \quad \cO_2= \{ {\bf a}_2, f_1{\bf a}_2\in {\bf\Sigma}_1, f^2_1{\bf a}_2 \in \Omega_2 , f^3_1{\bf a}_2,\dots, f^5_1{\bf a}_2 = \epsilon_1\}.$$
Combining (3.2--5) and (3.6--8), we have the matrix representations for $g$ and $f$:
$$ M_g = \pmatrix{2&1&0&0&0&0&0&1\cr 
                                 -1&-1&0&0&0&0&0&-1\cr 
                                  -1&-1&0&0&0&0&0&-1\cr
                                  0&0&1&0&0&0&0&0\cr 
                                  0&0&0&1&0&0&0&0\cr 
                                  0&0&0&0&1&0&0&0\cr 
                                  0&0&0&0&0&1&0&0\cr 
                                  0&0&0&0&0&0&1&0\cr}, \ 
        M_f= \pmatrix{2&1&0&0&0&0&0&1\cr 
                                 -1&-1&0&0&0&0&0&-1\cr 
                                  -1&-1&0&0&0&0&0&-1\cr
                                  0&0&1&0&0&0&0&0\cr 
                                  0&0&0&1&0&0&0&0\cr 
                                  -1&-1&0&0&1&0&0&-1\cr 
                                  0&-1&0&0&0&1&0&0\cr 
                                  0&0&0&0&0&0&1&0\cr}. $$
To show they have the same characteristic polynomial we are going to look at the matrices $M_g-xI$ and $M_f-xI$ where $I$ is the identity matrix and will show that after row and column operations to $M_f-xI$ we get the same matrix as $M_g-xI$.                                   
$$ M_f-xI= \pmatrix{2-x&1&0&0&0&0&0&1\cr 
                                 -1&-1-x&0&0&0&0&0&-1\cr 
                                  -1&-1&-x&0&0&0&0&-1\cr
                                  0&0&1&-x&0&0&0&0\cr 
                                  0&0&0&1&-x&0&0&0\cr 
                                  -1&-1&0&0&1&-x&0&-1\cr 
                                  0&-1&0&0&0&1&-x&0\cr 
                                  0&0&0&0&0&0&1&-x\cr}. $$
First we subtract the second row from the $6$-th row. For the general situation, we subtract rows for the lower generation chain from the rows for the corresponding part of the higher generation chain. Then we  add $6$-th column to the second column to remove the additional part and to obtain the matrix $M_g-xI$. For the general case, we add columns for the orbit collision part of the higher generation chain to the  corresponding lower generation chain to remove extra elements. It is clear that $M_f-xI$ and $M_g-xI$ have the same determinant. This gives us the desired result. 

\bigskip
\centerline{{\bf \S4. Non-degenerate Critical Triangle}}
\medskip
\noindent In this section we will determine the degree growth rate of $f$ with non-degenerate critical triangle.   As we noted at the beginning of \S2, it is equivalent to assume that 
$$ \beta_1 \ne 0, \quad {\rm and} \quad \beta_1\alpha_2-\alpha_1\beta_2 \ne 0. \eqno{(4.1)}$$
In particular, the curves $\Sigma_\gamma$, $\Sigma_\beta$ and $\Sigma_0$ are distinct, as well as $\{e_1,e_2,q\}$, the points of indeterminacy of $f^{-1}$.   Let us choose invertible linear maps $M_1$ and $M_2$ of $\P^2$ such that 
$$M_1{\bf\Sigma}_0=\Sigma_0,\ M_1{\bf \Sigma}_1=\Sigma_\beta,\ M_1{\bf \Sigma}_2=\Sigma_\gamma,$$
and 
$$M_2e_1=\epsilon_0,\ M_2e_2=\epsilon_1,\ M_2q=\epsilon_{2}.$$
It follows that $M_2\circ f_{\alpha,\beta}\circ M_1$ is a quadratic map with ${\bf \Sigma}_j\leftrightarrow e_j$ and so is equal to the map $J$.  Thus $f_{\alpha,\beta}$ is linearly conjugate to a mapping of the form $L\circ J$.  We will determine $\delta(\alpha,\beta)$ by finding the possibilities for $\#\cL^{c/o}$ and then applying Theorem 3.2.  When we treat $f_{\alpha,\beta}$ as a mapping $L\circ J$, we make the identifications 
$${\bf \Sigma}_0=\Sigma_0,\ \ {\bf \Sigma}_1=\Sigma_\beta,\ \ {\bf \Sigma}_2=\Sigma_\gamma,$$
$$\epsilon_0=p_\gamma,\ \ \epsilon_1=e_1,\ \ \epsilon_2=p_0,$$
and
$${\bf a}_0 = f({\bf\Sigma}_0-\cI(f))=e_1,\ \ {\bf a}_1 = f({\bf\Sigma}_1-\cI(f))=e_2,\ \ {\bf a}_2 = f({\bf\Sigma}_2-\cI(f))=q.$$

Thus  we have $f({\bf\Sigma}_0 - \cI)={\bf a}_0= \epsilon_1$, so the orbit $\cO_0= \{{\bf a}_0=\epsilon_1\}$ is singular and has length one.  There are two possibilities for the exceptional component ${\bf\Sigma}_1$; the first is that  ${\bf a}_1 \in {\bf\Sigma}_0 - \cI(f)$, which occurs when $\beta_2 \ne 0$.  (See Figure 4.1.)  The second possibility is ${\bf a}_1=\epsilon_2\in\cI$, which occurs when $\beta_2=0.$  (See Figure 4.2.)  An analysis of the possibilities for $\cO_{1}$ and $\cO_{2}$ will yield the candidates for $|\cO_{1}|$, $|\cO_{2}|$ and $\#\cL^{c/o}$, and thus give the possibilities for $\delta(\alpha,\beta)$.

\proclaim Theorem 4.1. If the critical triangle is non-degenerate and $\beta_2\ne0$, then $\delta_\star \leq \delta(\alpha,\beta) \leq \phi$. 

\noindent{\it Proof.}  Let $f_{Y}:Y\to Y$ be as in (1.2).  Since ${\bf a}_1=e_2 \ne \epsilon_i$ for $i = 0,1,2$, we have
$$ f_Y: \Sigma_1-\cI \to {\bf a}_1 \to [\beta_2:\alpha_2]_{E_1}\to [\beta_1\beta_2:\beta_1\alpha_2: \alpha_1\beta_2] \in \Sigma_B-\Sigma_0.\eqno{(4.2)} $$
If $f^2{\bf a}_1  =f_Y^2{\bf a}_1 = \epsilon_0$, then both lines $\Sigma_2$ and $\Sigma_B$ contain $\epsilon_1, \epsilon_0$. Since ${\bf a}_2 = \Sigma_B\cap \Sigma_1$ and $\epsilon_0 = \Sigma_2 \cap \Sigma_1$, we have ${\bf a}_2= \epsilon_0$. By the second statement of Lemma $1.1$, we see that the end points of both orbits $\cO_1$ and $\cO_2$ can not be $\epsilon_2$. It follows that we have at most two singular orbits including $\cO_0$. We have three cases. 

\epsfxsize=3.3in
\centerline{\epsfbox{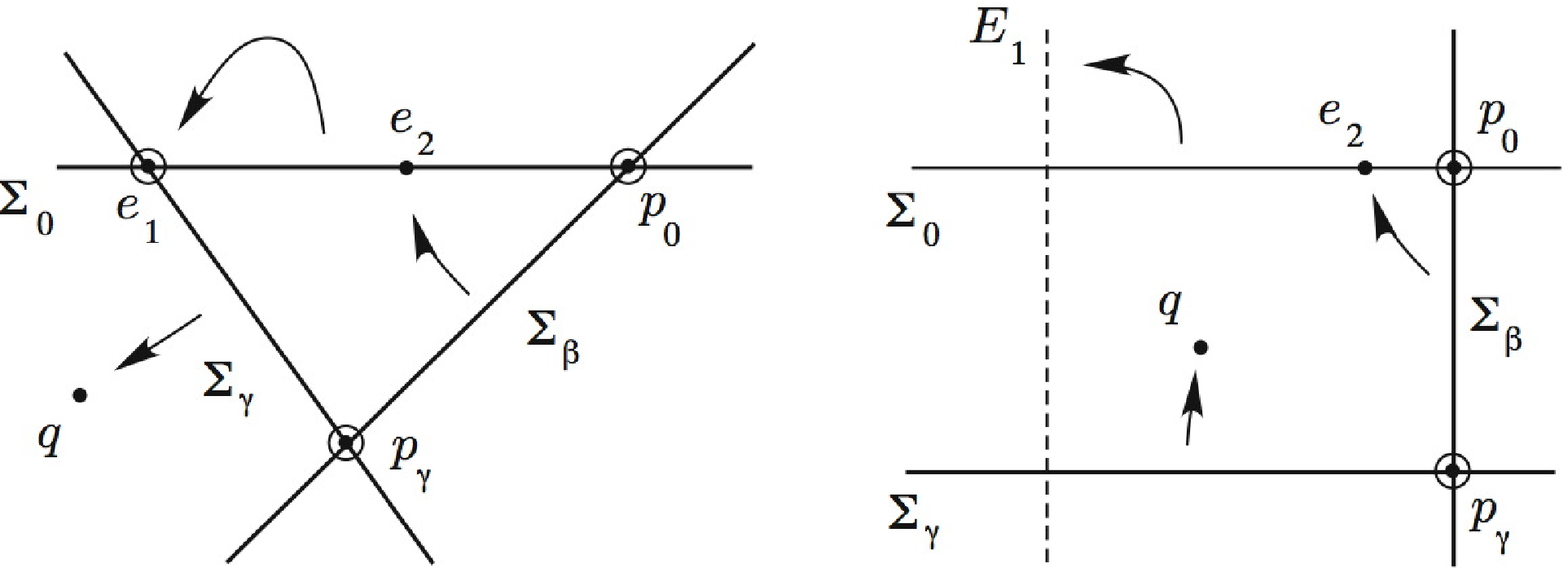}}
\centerline{Figure 4.1.  Nondegenerate critical triangle: case $\beta_2\ne0$.}
\medskip

The first case is where neither $\cO_1$ nor $\cO_2$ is singular. In this case the orbit list structure is $\#\cL^c = \emptyset, \#\cL^o=\{1\}$. By Theoerm 3.2,  $\delta(\alpha,\beta)$ is the largest real root of the polynomial 
$$\chi(x)=(x-2)x+(x-1) = x^2-x-1\eqno{(4.3)}$$ 
and is thus equal to $\phi$.

In the second case both $\cO_0$ and $\cO_1$ are singular. In this case the orbit $\cO_2$ can not be singular and therefore $f^2{\bf a}_1 \ne \epsilon_0$. By the equation $(4.2)$ with above argument, we have $n_1=|\cO_1|\geq 4$ and $\cO_1= \{{\bf a}_1, \dots, \epsilon_0\}$. It follows that $\#\cL^o =\emptyset, \#\cL^c = \{1,n_1\}$.  The dynamic degree $\delta(\alpha,\beta)$ is the largest root of the polynomial 
$$\chi(x)=(x-2)(x^{1+n_1}-1)+(x-1)(x+x^{n_1}+2)= x^{n_1}(x^2-x-1)+x^2.\eqno{(4.4)}$$ When $n_1=4$, the characteristic polynomial is given by $x^6-x^5-x^4+2 = x^2(x-1)(x^3-x-1)$.   Thus $\delta=\delta_{\star}$ in this case.  Let us observe that the Comparison Principle [BK, Theorem 5.1] concerns the modulus of the largest zero of the characteristic polynomial of $f^*$.  In \S3 we showed that the characteristic polynomials are the same in the elementary and the non-elementary cases.  Thus we may apply the Comparison Principle to conclude that $\delta(\alpha,\beta) \geq \delta_\star$ if $n_{1}\ge4$. 
\bigskip
\epsfxsize=3.5in
\centerline{\epsfbox{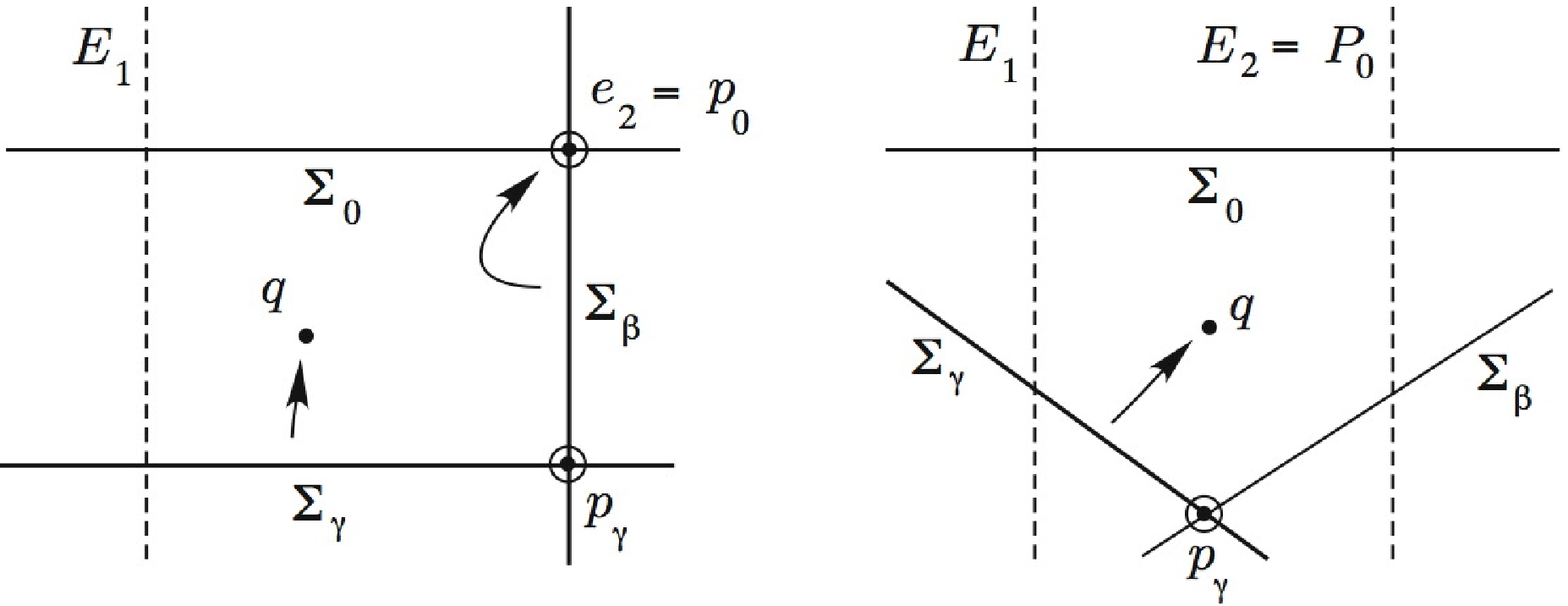}}
\centerline{Figure 4.2.  Nondegenerate critical triangle: case $\beta_2=0$.}
\medskip

The last case is where both $\cO_0$ and $\cO_2$ are singular. We have $n_2= |\cO_2|\geq 1$ and $\cO_2= \{{\bf a}_2, \dots, \epsilon_0\}.$ Therefore the orbit list structure is $\#\cL^c=\emptyset, \#\cL^o= \{n_2, 1\}$.  By Theorem 3.2, the dynamic degree $\delta(\alpha,\beta)$ is the largest root of the polynomial 
$$\chi(x)=(x-2)x^{1+n_2}+(x-1)(x+x^{n_2}+1)= x^{n_2}(x^2-x-1)+x^2-1.\eqno{(4.5)}$$ 
If $n_2 =1$, we have $\chi(x)= x^3-x-1$. 
\proclaim Theorem 4.2. Assume that the critical triangle is non-degenerate. If $\beta_2=0$ and $n_2=|\cO_2| \geq 8$, then $1< \delta(\alpha,\beta) \leq \delta_\star .$ If  $\beta_2=0$ and $n_2=|\cO_2| \leq 7$, then $\delta(\alpha,\beta)=1.$

\noindent{\it Proof.}  If $\beta_2=0$, we have ${\bf a}_1= \epsilon_2$ and therefore we have 
$$\cO_0=\{{\bf a}_0=\epsilon_1\}, \quad {\rm and}\quad \cO_1=\{{\bf a}_1=\epsilon_2\}.$$ If the orbit $\cO_2$ is non-singular, we have the orbit list structure $\#\cL^o =\{ 1,1\}, \#\cL^c=\emptyset$.   By Theorem 3.2, the degree growth rate $\delta(\alpha,\beta)$ is the largest root of the polynomial 
$$\chi(x) = (x-2)x^2+(x-1)(x+x+1)= x^3-x-1. \eqno{(4.6)}$$
If the orbit $\cO_2$ is singular, the end point of the orbit has to be the remaining point of indeterminacy, $\epsilon_0$. Thus we have $n_2= |\cO_2|\geq 1$ and $\cO_2= \{{\bf a}_2, \dots, \epsilon_0\}$. It follows that the orbit list structure $\#\cL^c =\{ 1,1, n_2\}, \#\cL^o=\emptyset$. Using the Lemma 2 and Proposition 7 in \S3, the dynamic degree is the largest root of the polynomial 
$$\eqalign{\chi(x) =& (x-2)(x^{2+n_2}-1)+(x-1)(2x^{1+n_2}+x^2+x^{n_2}+2x+3)\cr=&x^{n_2}( x^3-x-1)+x^3+x^2-1.\cr}\eqno{(4.7)}$$
It follows that $1\leq\delta(\alpha,\beta) \leq \delta_\star$. For $n_2=7$, we have $\chi(x)=(x^2-1)(x^3-1)(x^5-1)$ and so the $\delta(\alpha,\beta)=1$. For $n_2=8$ we have $\chi(x)= (x-1)(x^{10}+x^9-x^7-x^6-x^5-x^4-x^3+x+1)$ and $\chi'(1) <0$ and therefore the largest real root is strictly bigger than $1$.  It follows, then, from the comparison principle ([BK, Theorem 5.1]) that $\delta(\alpha,\beta)>1$ if $n_{2}\ge8$.

Let us note that when the orbit of $q$ lands on $p$, and we blow up the orbit of $q$, then we have removed the last exceptional curves for $f$ and $f^{-1}$.  Thus we have:

\proclaim Proposition 4.3.  If $(\alpha,\beta)\in V_n$, then the induced map $f_X:X\to X$ is biholomorphic.

Figure 4.3 shows the arrangement of the exceptional varieties in $X$ in the case where the orbit of $q$ does not enter $\Sigma_\beta$.

\epsfxsize=1.5in
\centerline{\epsfbox{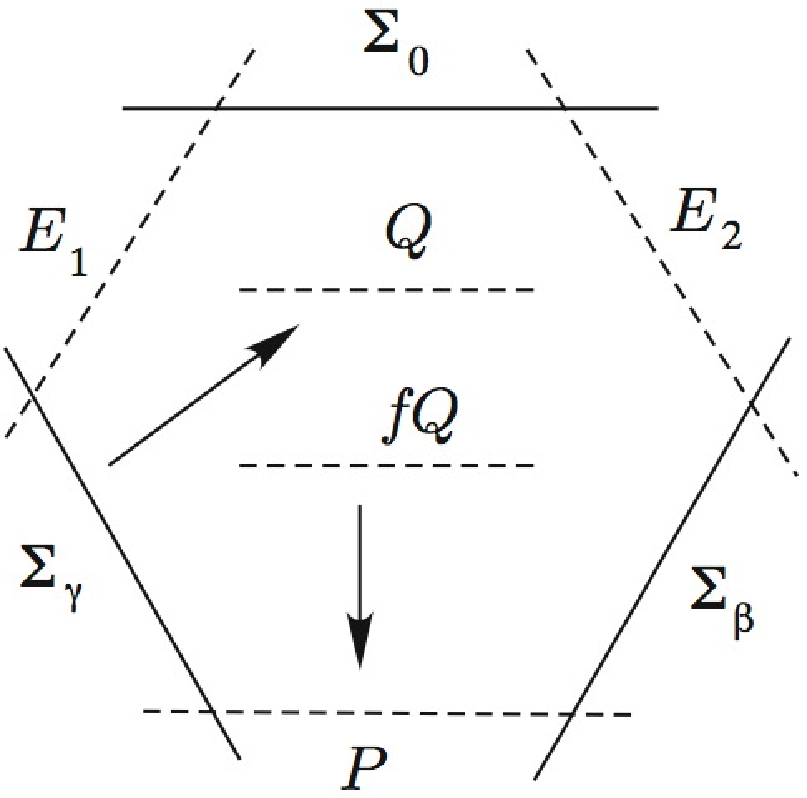}}
\centerline{Figure 4.3.  Nondegenerate critical triangle; elementary case $(\alpha,\beta)\in V_{n}$.}

\bigskip
\centerline{\bf \S5. Periodic Mappings   }
\medskip
\noindent  Here we determine the precise degree growth rate when $|\cO_{2}|\le7$.  In particular, we show that the degree grows quadratically when $|\cO_{2}|=7$, and we show that $f$ is periodic when $|\cO_{2}|\le 6$.  We do this by showing that $f^{*}$ is periodic in this case, and then we show that the periodicity of $f^{*}$ implies the periodicity of $f$.

Notice that if $|\cO_2|=n$, then $f^n(\Sigma_\gamma)=f^{n-1}(q)=p$, and therefore $(\alpha,\beta)\in V_{n-1}$.  To show the periodicity of $f_X^*$ it suffices to show that all roots of (4.7) with $n\le6$ are roots of unity and are simple.  For $n\le 6$ we list the characteristic polynomials, together with the smallest polynomials of the form $x^m-1$ that they divide:
$$\eqalign{
V_0\ (n=1):&\ \ (x-1)(x+1)(x^2+x+1)|(x^6-1)\cr
V_1\ (n=2):&\ \  (x-1)(x^4+x^3+x^2+x+1)|(x^5-1)\cr
V_2\ (n=3):&\ \  (x-1)(x+1)(x^4+1)|(x^8-1)\cr
V_3\ (n=4):&\ \  (x-1)(x^2+x+1)(x^4-x^2+1) |(x^{12}-1)\cr
V_4\ (n=5):&\ \ (x-1)(x+1)(x^6-x^3+1)|(x^{18}-1)\cr
V_5\ (n=6):&\ \  (x-1)(x^8+x^7-x^5-x^4-x^3+x+1)|(x^{30}-1)\cr}$$
Thus we have:
\proclaim Lemma 5.1.  Assume that the critical triangle is non-degenerate.  If $\beta_2=0$ and $n=|\cO_2|\le6$, then $f^*_X$ is periodic, with period $\kappa_n$, where $\kappa_n=6, 5,8,12, 18, 30$ (respectively).

When $|\cO_2|=7$, the largest root of equation (4.7) is 1 and has multiplicity 3.  Whether $f$ is elementary or not, the matrix representation from \S3 has a $3\times 3$ Jordan block with eigenvalue 1.  This means that $f^*_X$ has quadratric growth, and we have:
\proclaim Lemma 5.2.  Assume that the critical triangle is non-degenerate.  If $\beta_2=0$ and $|\cO_2|=7$, then $f^*_X$ has quadratic growth.

Notice that $|\cO_2|=1$ if and only if $q=p_\gamma$, which means that the parameters in $V_0$ satisfy  $\alpha_1\beta_0-\alpha_0\beta_1=-\alpha_2\beta_0=\alpha_1\alpha_2$.  With these conditions on $\alpha$ and $\beta$, it is not hard to check that the map $f$ is indeed periodic with period 6.  We could also see this by observing that $f$ has a period 6 cycle $\Sigma_\beta\mapsto e_2\mapsto \Sigma_0\mapsto e_1\mapsto\Sigma_\gamma\mapsto p_\gamma\mapsto\Sigma_\beta$.

\proclaim Theorem 5.3.  Assume that the critical triangle is non-degenerate.  if $\beta_2=0$ and $|\cO_2|\le 6$, then $f$ is periodic with period $\kappa_n$.

To prove Theorem 5.3, we use the following lemma:

\proclaim Lemma 5.4. If $f: \P^2 \to \P^2$ is a linear map with five invariant lines which are in general position, then $f$ is the identity.

\noindent{\it Proof.}  Let $l_i$, $i=0,1,2,3,4$, denote the lines fixed by $f$.  Since they are in general position, we may assume that they have the form
$ \Sigma_i=\{ x_i=0\} $ for $i=0,1,2$, 
$\Sigma_* = \{ x_0+w_1x_1+w_2x_2 = 0,\,\, w_1 \ne 0\}$, and 
$\Sigma_{**} = \{ x_0+v_1x_1+v_2x_2 = 0,\,\, v_2 \ne 0\}$.  A computation then shows that $f$ must be the identity.
\medskip
\noindent{\it Proof of Theorem 5.3.}  It suffices to show that $f^{\kappa_n}$ has at least five invariant lines for $n=2,\dots,6$.  Consider the basis elements $E_1$, $E_2$, $\cF_q$, and $\cF_{p_\gamma}$.  Since $(f^*_X)^{\kappa_n}$ is the identity, it fixes these basis elements.  Thus $f^{\kappa_n}$ fixes the base points in ${\bf P}^2$.  Since $f^{\kappa_n}$ is linear, it leaves invariant every line through two of these base points. 

\bigskip\noindent\centerline{\bf \S6.  Parameter Regions}
\medskip
\noindent  There is a natural group action on parameter space.  Namely, for $(\lambda,c,\mu)\in\C_{*}\times\C_{*}\times\C$ we have actions
$$(\alpha,\beta)\mapsto(\lambda\alpha,\lambda\beta)\eqno(6.1)$$
$$(\alpha,\beta)\mapsto (\alpha_{0},c\alpha_{1},c\alpha_{2},c\beta_{0},c^{2}\beta_{1},c^{2}\beta_{2})\eqno(6.2)$$
$$\eqalign{ (\alpha,\beta)\mapsto & (\alpha_{0}+\mu(\alpha_{1}+\alpha_{2})-\mu(\beta_{0}+\mu(\beta_{1}+\beta_{2})),  \cr & \alpha_{1}-\mu\beta_{1},\alpha_{2}-\mu\beta_{2},\beta_{0}+\mu(\beta_{1}+\beta_{2}),\beta_{1},\beta_{2}).\cr}\eqno(6.3)$$
The first action corresponds to the homogeneity of $f_{\alpha,\beta}$.  The other two are given by linear conjugacies of $f_{\alpha,\beta}$.  To see them, we write $f$ in affine coordinates, as in (0.2).  Action (6.2) is given by conjugating by the scaling map $(x_{1},x_{2})\mapsto(cx_{1},cx_{2})$, and (6.3) is given by conjugating by the translation $(x_{1},x_{2})\mapsto(x_{1}+\mu,x_{2}+\mu)$.   

Now consider maps of the form
$$f:(x,y)\mapsto (y,{y\over b+x+cy}),\ \ b\ne0.\eqno(6.4)$$
In this case we have $\alpha=(0,0,1)$, $\beta=(b,1,c)$ and $\gamma=(0,0,1)$.  Let $Y$ be as in (1.2), and let $f_{Y}:Y\to Y$ be the induced map.  Repeating the computation of (1.3), we see that
$$\Sigma_{\beta}\mapsto E_{2}\mapsto [c:0:1]_{e_{1}}\in E_{1}\mapsto [c:1:0]\in\Sigma_{\gamma}.\eqno(6.5)$$
We conclude that the sub-family (6.4) is critically finite the following sense that all exceptional curves have finite orbits:
\proclaim Proposition 6.1.  If $f$ be as in (6.4), then $q=(0,0)$ is a fixed point, and the exceptional curves are mapped to $q$.  In particular, $f_{Y}$ is 1-regular.

\noindent{\it Proof.}  If $c=0$, then exceptional locus is $\Sigma_{\gamma}$; if $c\ne0$, then both $\Sigma_{\beta}$ and $\Sigma_{\gamma}$ are exceptional.  We see from (6.5) that in either case the exceptional curves are mapped to the fixed point.
\medskip
The variety $V_{n}\subset\{\beta_2=0\}$ corresponds to a dynamical property: an exceptional line is  mapped to a point of indeterminacy.  Thus $V_{n}$ is invariant under the actions (6.1--3).  For $(\alpha,\beta)\in V_{n}$, we have $\beta_{2}=0$, and we may apply (6.3) to obtain $\alpha_{1}=0$.  Since by (0.4) we must have $\alpha_{2}\ne0$ and $\beta_{1}\ne0$, we apply (6.1) and (6.2) to obtain $\alpha_{2}=\beta_{1}=1$.  Thus each orbit within $V_{n}$ is represented by a map which may be written in affine coordinates as
$$(x,y)\mapsto(y,{a+y\over b+x}).\eqno(6.6)$$
If $f$ is of the form (6.6), then $f^{-1}$ is conjugate via the involution $\sigma:x\leftrightarrow y$ and a transformation (6.3) to the map 
$$(x_,y)\mapsto(y,{a-b+y\over -b+x}).\eqno(6.7)$$  
Such a mapping is conjugate to its inverse if $b=0$.  

Now we suppose that  $f$ is given by (6.6).  Thus $q=(-a,0)$ and $p=(-b,-a)$, and $V_n$ is defined by the condition $f^nq=p$.
The coefficients of the equations defining $V_{n}$ are positive integers, and $V_n$ is preserved under complex conjugation.  An inspection of the equations defining $V_{n}$ produces the first few:
\smallskip
$V_{0}$:  the orbit under (6.1--3) of $(a,b)=(0,0)$
\smallskip
$V_{1}$:  the orbit of $(a,b)=(1,0)$
\smallskip
$V_{2}$:  the orbits of $(a,b)=((1 + i)/2, i)$ and its conjugate.
\smallskip
$V_{3}$:  the orbits of $(a,b)\in\{(2 +i -\sqrt{3})/2, i),(2 +i +\sqrt{3})/2, i)\}$ and their conjugates.
\medskip
\noindent We solve for $V_4$, $V_5$ and $V_6$ by using the resultant polynomials of the defining equations, and we find:
\smallskip
$V_4$: the orbits of $(a,b)=(0.8711+ 0.7309i,1.4619i)$,
$(0.6974 + 0.2538i, 0.5077i)$, 
$(-0.06857+ 0.3889i,  0.7778i)$,
and their conjugates.  The exact values  are roots of $1 - 3\ a + 9a^2 - 24a^3 + 36a^4 - 27a^5 + 9a^6$ and $1 + 6b^2 + 9b^4 + 3b^6$.
\smallskip
$V_5$: the orbits of $(a,b)=(3.7007+ 1.2024, 2.4048i)$,
$(1.0353+ 0.3364i,0.6728i)$,
$( 0.4465+ 0.6146i,1.2293i)$,
$( -0.1826+ 0.2513i, 0.5027i),$
or their conjugates.  The exact values are roots of $1 + 3a^2 - 20a^3 + 49a^4 - 60a^5 + 37a^6 - 10\ a^7 + a^8$ and $1 + 7b^2 + 14b^4 + 8 b^6 + b^8$.
\smallskip
$V_6$:  The defining equations for $V_6$ are divisible by $b^2$, so all points of the form $(a,0)$, $a\ne0,1$, belong to $V_{6}$.  By (6.7), these parameters correspond to maps which are conjugate to their inverses.   In addition, $V_{6}$ contains the orbits of
$$a=(3\pm\sqrt5+2b)/4,\ \ b=i\sqrt{(5\pm\sqrt 5)/2} $$
and their conjugates.  

By Theorem  2, mappings in $V_6$ have quadratic degree growth, and by [G] such mappings have invariant fibrations by elliptic curves.  Let us show how our approach yields these invariant fibrations.  
\medskip
\centerline{\epsfysize1.6in\epsfbox{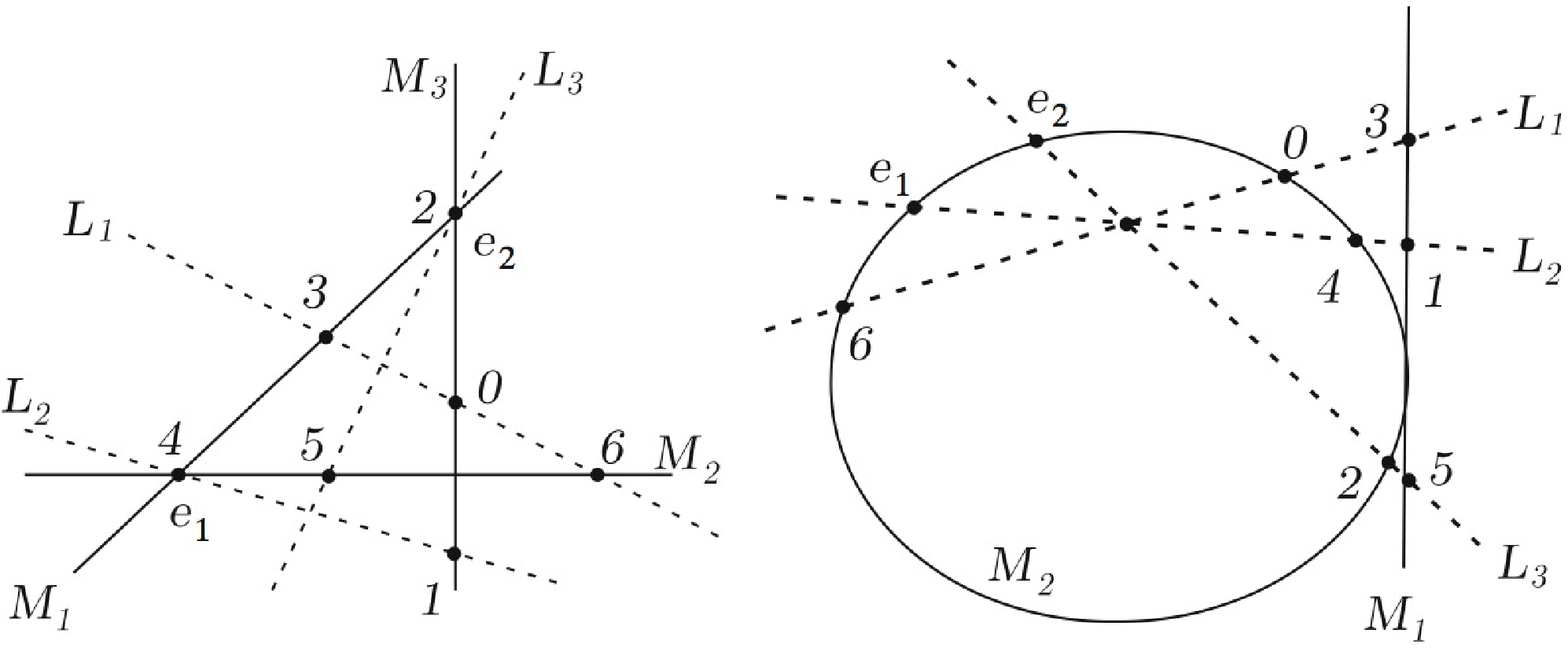}}
\centerline{Figure 6.1.  Points $f^jq$ = `$j$', $0\le j\le 6$, for $V_6$.  Case $b=0$ on left; $b\ne0$ on right.}
\bigskip

Let us first consider parameters $(a,0)$.  In this case, the fibration was obtained classically in [L] and [KoL].    In the space $Y$ of (1.2), the $f$-orbit $\{q_{j}=f^{j}q:j=0,1,\dots,6\}$ is:
$$\eqalign{&q_{0}=(-a,0)_{{\bf C}^2}=[1:-a:0],\ q_{1}=(0,-1)_{{\bf C}^2}=[1:0:-1],\ q_{2}=[0:0:1]=e_2,\cr
& q_{3}=[0:1:-1],\ q_{4}=[1:0:-1]_{e_1},\ q_{5}=(-1,0),\ q_{6}=(0,-a)=p,\cr}$$
as is shown in Figure 6.1.  Here we use `$j$' to denote `$q_j$'.  The construction of $X$ is shown in Figure~6.2, where `$f^jQ$' denotes the blowup fiber over $q_j$.  In contrast,  the case corresponding to $(a,b)\in V_6$, $b\ne0$ corresponds to Figure 4.3.  Consulting Figure 6.2, we see that the cohomology class $3H_X-E_1-E_2-Q_2-Q_4-\sum Q_j$ is fixed under $f^*$.  We will find polynomials which correspond as closely as possible to this class.  These will be cubics which vanish on $e_i$ and $q_j$.  
Looking for lines that contain as many of the $q_j$ as possible, we see $L_{1}=\{x+y+a=0\}$ contains 0,3,6.  Mapping forward by $f$, we have 
$$L_{1}\mapsto L_{2}=\{y+1=0\}\mapsto L_{3}=\{x+1=0\}\mapsto L_{1}.$$
In addition, the points $q_{j}$, $j=2,3,4$ are contained in the line at infinity $M_{1}=\Sigma_{0}$.  This maps forward as:
$$M_{1}\mapsto e_{1}\mapsto M_{2}=\{y=0\}\mapsto M_{3}=\{x=0\}\mapsto e_{2}\mapsto M_{1}.$$
The cubic $c_{1}=(x+y+at)(x+t)(y+t)$ defines $L_{1}+L_{2}+L_{3}$ in ${\bf P}^{2}$, and $c_{2}=xyt$ defines $M_{1}+M_{2}+M_{3}$.  Setting $t=1$ and taking the quotient, we find the classical invariant $h(x,y)=c_{1}/c_{2}$.

\centerline{\epsfysize1.6in\epsfbox{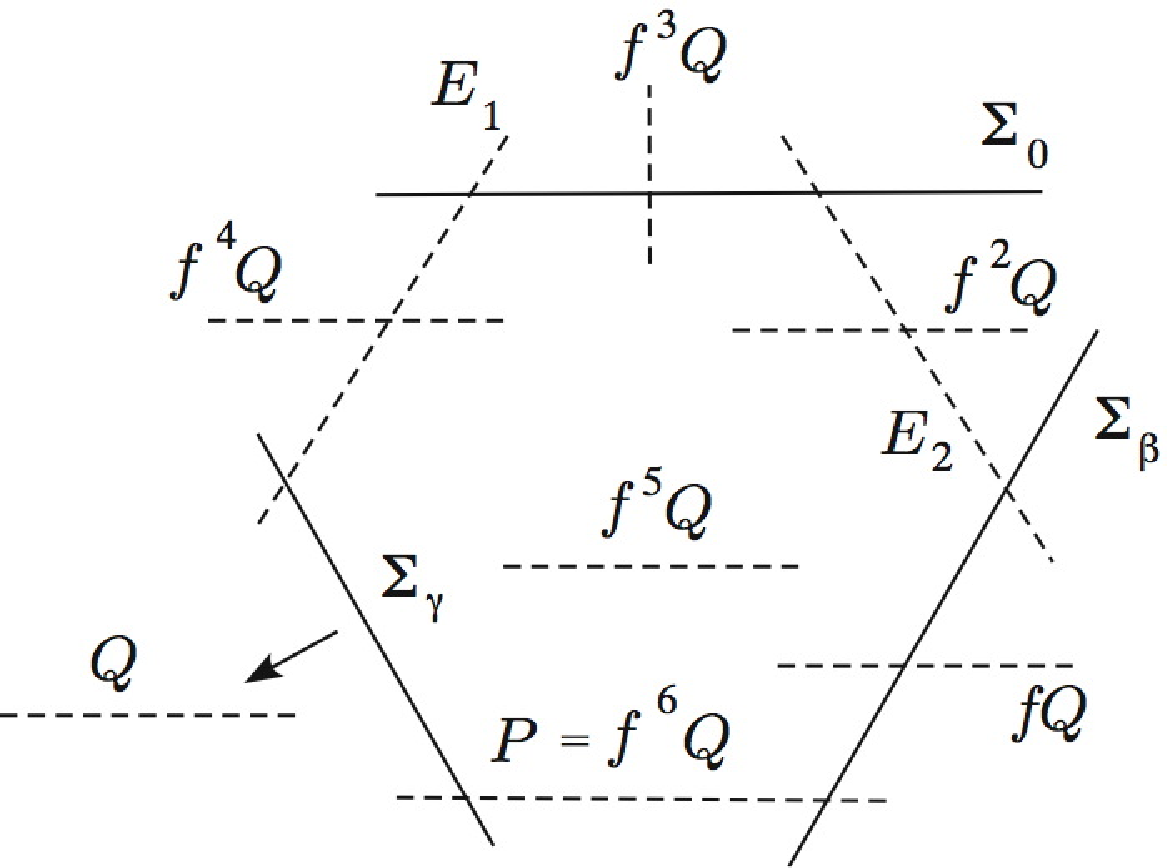}}
\centerline{Figure 6.2.  Space $X$ for $V_6$, $b=0$.}
\bigskip

Now we consider the other four parameters $(a,b)$ in $V_{6}$.  Inspecting the defining equations of $V_{6}$, we find that $a$ and $b$ satisfy $-2a+a^{2}+b-ab=0$ and $-b^{2}-1+b-2a=0$.  Using these relations, we see that the $f$-orbit of $q$ is:
$$\eqalign{q_{0}=(-a,0), \ q_{1}&=(0,1-a),\ q_{2}=(1-a,1/b), \ q_{3}=(1/b,a(1+ab)/(ab-b^2)),\cr
 q_{4}&=(a(1+ab)/(ab-b^2),1-a),\ q_{5}=(1-a,-b),\ q_{6}=(-b,-a).\cr}$$
 Looking again at the points $q_{j}$, $j=0,3,6$, we see that they are contained in a line $L_{1}=\{x+(1-{b\over a})y + a=0\}$.  Mapping $L_{1}$ forward under $f$, we find:
$$L_{1}\mapsto L_{2}=\{y+a-1=0\}\mapsto L_{3}=\{x+a-1=0\}\mapsto L_{1}.$$
We multiply these linear functions together to obtain a cubic $c_{1}$ which defines  $\sum L_{i}$.
We see, too, that the points $q_{j}$, $j=1,3,5$ are contained in the line $M_{1}=\{(a-b-1)x+(a-1)y+(a-1)^{2}=0\}$.  Mapping forward, we find:
$$M_{1}\mapsto M_{2}=\{(a-1)xy +(b^{2}+1)y+(a-b-1)x+(a-b)=0\}\mapsto M_{1}.$$
Multiplying the defining functions, we obtain a cubic $c_{2}$ which defines $M_{1}+M_{2}$.  Now we define
$k(x,y)=c_{1}/c_{2}$.  And inspection shows that $k\circ f=\omega k$, where $\omega$ is a 5th root of unity.  Thus $f$ is a period 5 mapping of the set of cubics $\{k={\rm const}\}$ to itself.

\bigskip\centerline{\bf Appendix:  Explanation of the Computer Graphics}

\medskip\noindent  It is useful to have visual representations for rational mappings.  A number of interesting computer graphic representations of the behavior of rational mappings of the real plane have been given in various works by Bischi, Gardini and Mira; we cite [BGM] as an example.  The pictures here have a somewhat different origin and are made following a scheme used earlier by one of the authors and Jeff Diller (see [BD1,3]).  They are motivated by the theory of dynamics of complex surface maps.  Let $f$ be a birational map of a K\"ahler surface.  If $\delta(f)>1$, then there are positive, closed, (1,1)-currents $T^{\pm}$ such that $f^{*}T^{+}=\delta(f)T^{+}$ and $f^{*}T^{-}=\delta(f)^{-1}T^{-}$ (see Diller-Favre [DF]).  These currents have the additional property that for any complex curve $\Gamma$ there is a number $c>0$ such that
$$cT^{+} = \lim_{n\to\infty}{1\over\delta^{n}}f^{n*}[\Gamma],\eqno(A.1)$$
and similarly for $T^{-}$.  By work of Dujardin [D1] these currents have the structure of a generalized lamination.  We let  ${\cal L}^{s/u}$ denote the generalized laminations corresponding to $T^{\pm}$.  It was shown in [BD2] that the wedge product $T^{+}\wedge T^{-}$ defines an invariant measure in many cases, and Dujardin [D2] showed that this invariant measure may be found by taking the ``geometric intersection'' of the measured laminations ${\cal L}^{s}$ and ${\cal L}^{u}$.

When one of our mappings $f$ has real coefficients, it defines a birational map of the real plane, and we can hope that there might be real analogues for the results of the theory of complex surfaces.  This was proved to be the case for certain maps in [BD1,3] but is not known to hold for the maps studied in the present paper.

Figure 0.1 was drawn as follows.  We work in the affine coordinate chart $(x,y)$ on ${\bf R}^2$ given by $x_0=1$, $x=x_1/x_0=x_1$, $y=x_2/x_0=x_2$.  We start with a long segment $L\subset{\bf R}^{2}$ and map it forward several times.  The resulting curve is colored black and ``represents'' ${\cal L}^{u}$.  After the first few iterates, the computer picture seems to ``stabilize,'' and further iteration serves to ``fill out''  the lamination. The appearance of the computer picture obtained in this manner is independent of the choice of initial line $L$.  To represent ${\cal L}^{s}$, we repeat this procedure for $f^{-1}$ and color the resulting picture gray.  In Figure 0.1 we present ${\cal L}^s$ in gray in the left hand frame.  Then we present ${\cal L}^s$ and ${\cal L}^u$ together in the right hand frame in order to show the set where they intersect.

\centerline{\epsfxsize2.4in\epsfbox{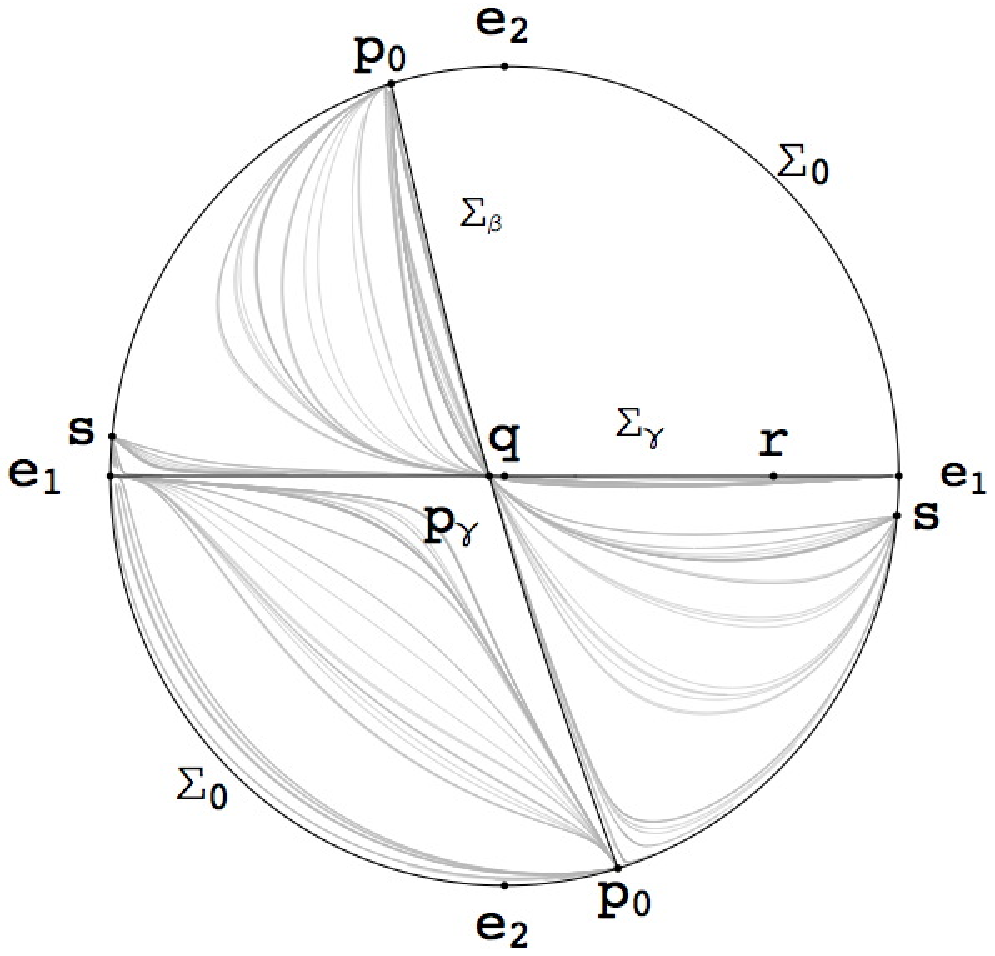}\hfil \epsfxsize2.4in\epsfbox{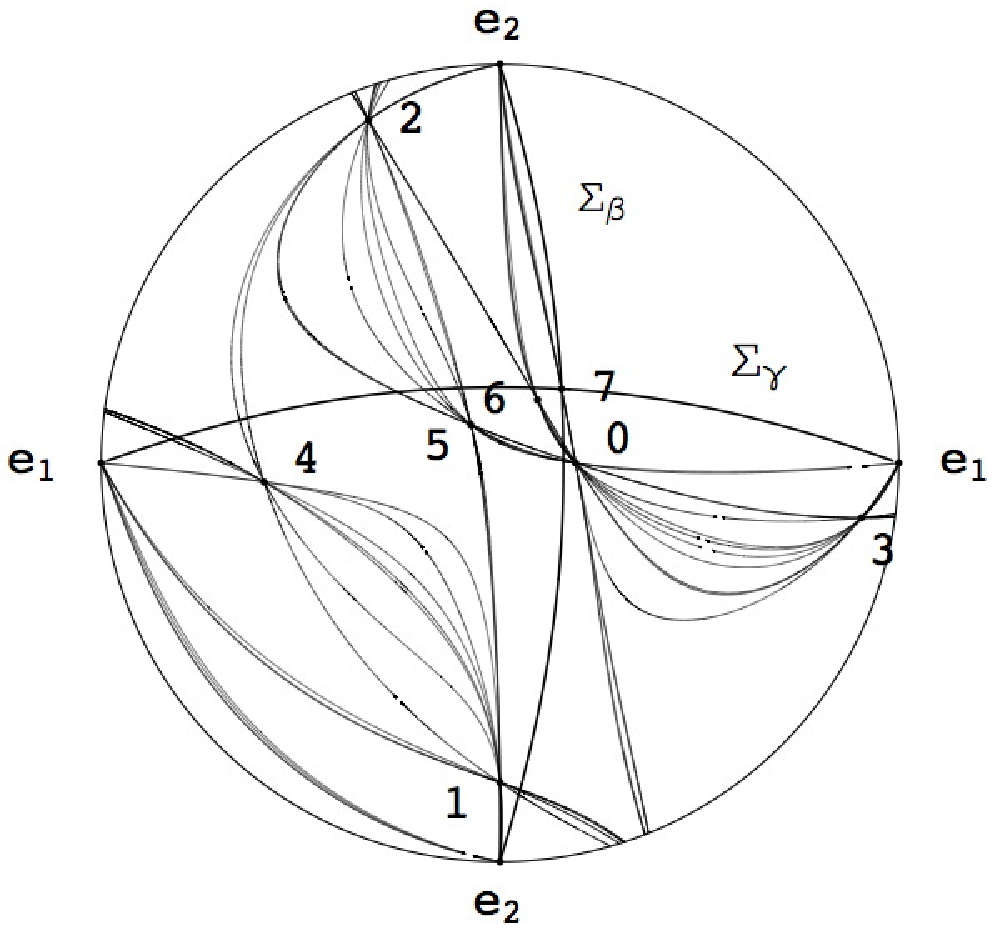}}
\centerline{Figure A.1.  Explanation of Figure 0.1 (left); a mapping from $V_{7}$ (right).}
\medskip

We also want the graphic to have the appearance of a subset of ${\bf P}^{2}$,  so we rescale the distance to the origin.  The resulting ``disk'' is a compactification of ${\bf R}^{2}$.  In fact, this is real projective space, since antipodal points of the circle are identified.  The circle forming the boundary of this disk is the line at infinity $\Sigma_0$.   

Figure 0.1 is obtained using the map of the form (6.4):
$$(x,y)\mapsto (y,{y\over .1+x+.3y}).$$
By Proposition 6.1, $f$ is critically finite, so $\delta(f)=\phi$ by Theorem 4.1.  On the left half of 
Figure A.1, we have re-drawn ${\cal L}^{s}$, together with the points of indeterminacy of $f$ and $f^{-1}$.   Pictured, for instance, are $e_1$, $e_2$, $p_0=[0:-.3:1]$, $p_{\gamma}=(-.1,0)$, and $q=(0,0)$. The exceptional curves are lines connecting certain pairs of these points and may be found easily using Figure 1.1 as a guide.  As we expect,  ${\cal L}^s$ is ``bunched'' at the points of indeterminacy of $f$, i.e., $p_0$, $e_1$, and $p_\gamma$.  Let us track the backward orbits of these points.  First, $p_0=f^{-1}p_0$ is fixed under $f^{-1}$, and $f^{-1}p_\gamma=e_1$.  Now let $Y$ and $f_Y$ be as in (1.2).  Repeating the calculations at equation (1.3), we see that $f_Y^{-1}$ takes $p_\gamma$ to the fiber point $[1:0:-.1]_{E_1}$ over $e_1$.  Then this fiber point is mapped under $f^{-1}$ to the point $s=[0:1.03:-.1]\in\Sigma_0$.  The next preimage is $f^{-1}s=p_0$, so $f^{-1}$ is critically finite in the sense that the exceptional curves all have finite orbits.  This explains why ${\cal L}^s$ is ``bunched'' at only four points. 

To explain the points where ${\cal L}^u$ is ``bunched,'' we have plotted the point $r:=f^3\Sigma_\beta=(10/3,0)$ from (6.5).  If we superimpose the picture of ${\cal L}^u$ on the left panel of Figure A.1, we find that ${\cal L}^{u}$ is ``bunched''  exactly on the set $e_1$, $e_2$, $q$, and $r$.
The ``eye'' which appears in the first quadrant is due to an attracting fixed point.

The right hand side of Figure A.1 is obtained using the map
$$(x,y)\mapsto(y,{-.499497+y\over -.415761+x}),$$
which corresponds to a real parameter $(a,b)\in V_7$.  By  ``$j$'', $j=0,\dots,7$, we denote the point  $f^jq$.  Thus ``7'' is the point of indeterminacy $p=f^7q$.  We let $\pi:X\to{\bf P}^2$ be the manifold obtained by blowing up $e_1$, $e_2$, and ``$j$'' for $j=0,\dots,7$.  The lamina of ${\cal L}^u$ are then separated in $X$, and the apparent intersections may be viewed as artifacts  of the projection $\pi$.

\bigskip\noindent\centerline{\bf References}
\medskip

\item{[BD1]} E. Bedford and J. Diller,  Real and complex dynamics of a family of birational maps of the plane:  the golden mean subshift, American J. of Math., to appear.

\item{[BD2]} E. Bedford and J. Diller,  Energy and invariant measures for birational surface maps, Duke Math.\ J., to appear.

\item{[BD3]} E. Bedford and J. Diller,  Dynamics of a two parameter family of plane birational mappings: Maximal entropy.  arxiv.org/math.DS/0505062

\item{[BK]} E. Bedford and KH. Kim, On the degree growth of birational mappings in higher dimension,  J. Geom.\  Anal.\ 14 (2004), 567-596. 

\item{[BDGPS]} M.J. Bertin, A. Decomps-Guilloux, M. Grandet-Hugot, M. Pathiaux-Delefosse, and J.P. Schreiber, {\sl Pisot and Salem Numbers}, Birkh\"auser Verlag (1992)

\item{[BGM]}  G-I. Bischi, L. Gardini, and C. Mira,  Plane maps with denominator.  I.  Some generic properties, International J. of Bifurcation and Chaos, 9 (1999), 119--153.

\item{[C]} S. Cantat, Dynamique des automorphismes des surfaces projectives complexes,  C.R.\ Acad.\ Sci.\ Paris, t.\ 328, p.\ 901--906, 1999.

\item{[CL]} M. Cs\"ornyei and M. Laczkovich, Some periodic and non-periodic recursions, Monatshefte f\"ur Mathematik 132 (2001), 215-236. 

\item{[DF]}  J. Diller and C. Favre, Dynamics of bimeromorphic maps of
surfaces, Amer. J. of Math., 123 (2001), 1135--1169.

\item{[D1]} R. Dujardin,  Laminar currents in ${\bf P}^2$, Math.\ Ann., 325, 2003, 745--765.

\item{[D2]} R. Dujardin,  Laminar currents and birational dynamics.  arxiv.org/math.DS/0409557

\item{[FS]} J-E Forn\ae ss and N. Sibony,  Complex dynamics in higher dimension, II.  {\sl Modern Methods in Complex Analysis}, Ann. of Math. Studies, vol. 137, Princeton Univ. Press, 1995, pp. 135--182.

\item{[GBM]} L. Gardini, G.I. Bischi, and C. Mira, Invariant curves and focal points in a Lyness iterative process, Int. J. Bifurcation and Chaos 13 (2003), 1841-1852.

\item{[G]}  M. Gizatullin, Rational $G$-surfaces. (Russian) Izv. Akad. Nauk SSSR Ser. Mat. 44 (1980), no. 1, 110--144, 239. 

\item{[GKP]} R.L. Graham, D.E. Knuth, and O. Patashnik, {\sl Concrete Mathematics},  1989.

\item{[GL]}  E.A. Grove and G. Ladas,  {\sl Periodicities  in Nonlinear Difference Equations},  Kluwer Academic Publishers, 2005.

\item{[KoL]} V.I. Kocic and G. Ladas, {\sl Global Behaviour of Nonlinear Difference Equations of Higher Order with Applications}, Kluwer Academic Publishers 1993.

\item{[KLR]} V.I. Kocic, G. Ladas, and I.W. Rodrigues, On rational recursive sequences, J. Math. Anal. Appl 173 (1993), 127-157.

\item{[KuL]}  M. Kulenovic and G. Ladas,  {\sl Dynamics of Second Order Rational Difference Equations}, CRC Press, 2002.


\item{[KG]}  R.P. Kurshan and B. Gopinath, Recursively generated periodic sequences, Canad.\ J. Math. 26 (1974), 1356--1371.


\item{[L]} R.C. Lyness, Notes 1581,1847, and 2952, Math. Gazette {\bf 26} (1942), 62, {\bf 29} (1945), 231, and {\bf 45} (1961), 201.
\bigskip
\rightline{Indiana University}

\rightline{Bloomington, IN 47405}

\rightline{\tt bedford@indiana.edu}

\bigskip
\rightline{Syracuse University}

\rightline{Syracuse, NY 13244}

\rightline{\tt kkim26@syr.edu}

\end